\documentclass[11pt]{article}
\usepackage[top=1in,bottom=1in,left=1in,right=1in]{geometry}
\usepackage[english]{babel}
\usepackage{amsthm}
\usepackage{tabularx}
\usepackage{adjustbox}
\usepackage{graphicx}
\usepackage{epstopdf}
\usepackage{algorithmic}
\usepackage{wrapfig}
\usepackage[version=4]{mhchem}
\usepackage{amsmath}
\usepackage{amssymb}
\usepackage{hyperref}
\usepackage[T1]{fontenc}
\usepackage{soul,xcolor}
\usepackage{cases}
\usepackage{booktabs}
\usepackage{float}
\usepackage{bm}
\usepackage{svg}
\usepackage{setspace}
\usepackage{makecell}
\usepackage{longtable}
\usepackage{mathtools}
\usepackage{enumerate}
\usepackage{lineno}
\usepackage{pdfpages}
\usepackage{graphicx}
\usepackage{subcaption}
\usepackage{multirow}
\usepackage{array}        
\newcolumntype{Y}{>{\centering\arraybackslash}X}  
\usepackage{accents}
\usepackage{enumitem}
\usepackage{nomencl}
\makenomenclature
\renewcommand{\nomgroup}[1]{%
  \ifthenelse{\equal{#1}{A}}{\item[\textbf{\textit{Indices}}]}{}%
  \ifthenelse{\equal{#1}{B}}{\item[\textbf{\textit{Sets}}]}{}%
  \ifthenelse{\equal{#1}{C}}{\item[\textbf{\textit{Parameters}}]}{}%
  \ifthenelse{\equal{#1}{D}}{\item[\textbf{\textit{Variables}}]}{}%
}

\begin{document}

\title{Toward Decarbonization of Chemical Manufacturing: Joint Optimization of Unit Commitment and Microgrid Operations}

\author{
  Saba Ghasemi Naraghi$^1$ \and Richard Reed$^2$ \and Tylee Kareck$^1$ \and Paritosh Ramanan$^2$ \and Zheyu Jiang$^2$\\
}

\date{
    \normalsize $^1$School of Chemical Engineering, Oklahoma State University, Stillwater, OK 74078\\
    $^2$School of Industrial Engineering and Management, Oklahoma State University, Stillwater, OK 74078
}

\maketitle

\begin{abstract}
The electrification of chemical process heating is essential to industrial decarbonization and sustainable manufacturing of chemical products. Joint optimization of electrified chemical process heating units and electric power systems is needed to achieve decarbonized operation of both sectors. In this work, we introduce a joint optimization model that identifies the optimal unit commitment of power systems and optimal operation of electrified steam cracking microgrids for olefins production. A mixed-integer linear programming (MILP) model is developed to optimize the hourly operational plan of 26 ethylene plants and the main power system in Texas. We propose a two-stage solution method to solve the resulting large-scale MILP problem efficiently. In the first stage, we apply Benders decomposition with LP-relaxed subproblems to decouple microgrid operations from the main power system. In the second stage, we use the first-stage optimal solution as a warm-starting point for the MILP. This two-stage approach reduces the average solution time by 93.5\% compared to direct solution of the MILP. Results show that the largest overall greenhouse gas (GHG) emission reduction for both power systems and microgrids is achieved when the electrification level of steam cracking units is at 30\%. Above 30\% electrification level, a higher electrification level leads to higher overall GHG emissions and steady increase in operating costs, particularly on the microgrid side. Increasing renewables contributions in the electric power system helps debottleneck the electrification efforts and facilitate holistic decarbonization. We also remark that the optimal operational plan of electrified steam cracking microgrids also exhibit strong spatiotemporal patterns.

\textit{Keywords:} Microgrid operation, unit commitment problem, steam cracking, decarbonization and electrification, optimization.
\end{abstract}

\section{Background and Motivation}
The U.S. manufacturing sector accounts for approximately 20\% of the nation's total primary energy consumption and greenhouse gas (GHG) emissions \cite{USDOE2022}. The top two energy-use process industries in the U.S., chemical and refining, consume nearly half of the manufacturing sector's primary energy and emit half of its GHGs. A significant portion of this energy demand is for process heating, which uses thermal energy to convert feedstock into products. Currently, almost all process heating needs are met by direct combustion of fossil fuels or through the use of steam, which is generated directly or indirectly from the combustion of fossil fuels \cite{AgrawalSiirola2023}. However, as the U.S. energy landscape continues to shift toward clean and renewable electricity \cite{USDOE2022}, such as solar and wind, the energy can be directly harvested as electricity, making the electrification of process heating a viable and attractive solution to reducing carbon footprint and energy usage of chemical and refining processes \cite{mallapragada2023decarbonization}. 

As substantial industrial loads from electrified process heating units are integrated with a power grid that increasingly depends on variable renewable energy (VRE), the successful decarbonization of chemical and refining industries via the transition of process heating mechanisms depends on close coordination between chemical plants and power system stakeholders, including Independent System Operators (ISOs), to optimize their operations simultaneously, while ensuring that electrified process heating units can operate continuously at steady state despite the intermittency of VRE \cite{USEIA2025}.

To study and illustrate how such coordination can be optimally achieved and its economic and environmental impacts, in this work, we focus on a key application, ethylene production via steam cracking, and propose a joint optimization framework for day-ahead unit commitment of power systems and operational schedules of electrified steam cracking units. Ethylene is a platform chemical that serves as the key precursor for plastics, textiles, antifreeze, detergents, adhesives, rubber, and food packaging materials \cite{afpm}. The global annual production rate of ethylene exceeded 200 million metric ton (MMT) in 2022 and is expected to grow by more than 60\% by 2034 \cite{Precedence2024}. 

Ethylene is predominantly produced through steam cracking, an energy-intensive process that cracks hydrocarbon feedstocks such as ethane, propane, and naphtha at high temperatures in furnaces (known as cracker) \cite{Zimmermann2009}. The furnace heat is supplied by the combustion of fresh natural gas fuel or the methane fraction of the cracking byproduct, making steam cracking one of the most energy and emission-intensive chemical processes \cite{Zimmermann2009}. A promising pathway to decarbonize the steam cracking process is electrification \cite{USDOE2022}. For instance, in 2024, BASF, SABIC, and Linde launched the world's first demonstration plant for an electrified steam cracking furnace with 6 MW of renewable electrical energy input \cite{Sinn2024}, which can reduce \ce{CO2} emissions by at least 90\% \cite{sabicnews}. However, one critical challenge is that existing ethylene plants are designed to operate uninterruptedly in steady state, whereas VRE (solar and wind) contribution to U.S. electricity generation is projected to increase from 19\% in 2025 to more than 64\% by 2050 \cite{USEIA2025}. Furthermore, large-scale energy storage or complete plant redesigns would impose prohibitive economic and technical barriers to chemical companies, making immediate full-scale electrification impractical \cite{AgrawalSiirola2023}.

\begin{figure}[ht!]
    \centering    
    \includegraphics[width=\textwidth]{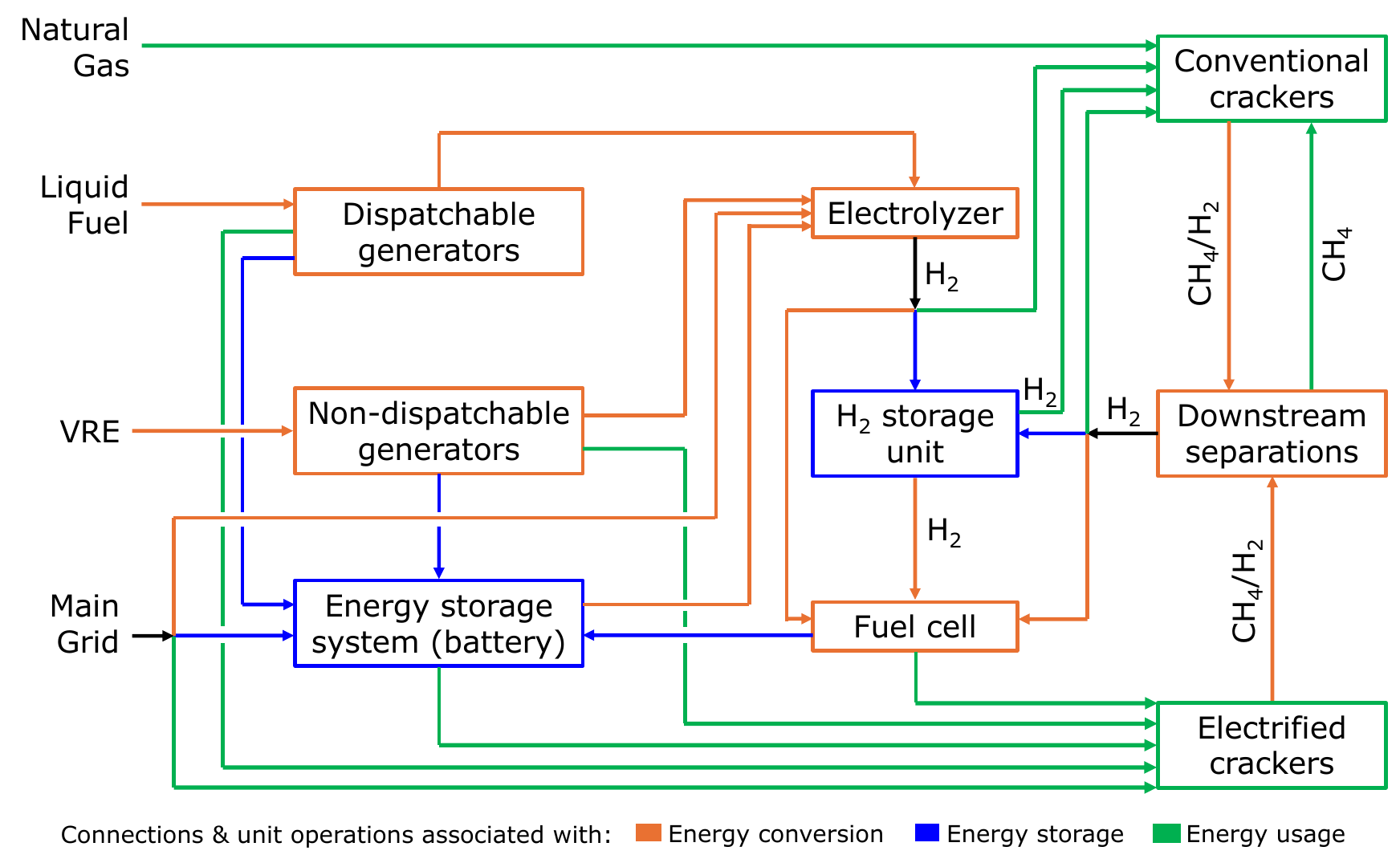}
    \vspace{-1em}
    \caption{Our envisioned framework for using electricity to supply process heat for steam cracking, which is adopted from our earlier work \cite{Naraghi2025}. Diverse energy sources supply heat for both electrified and conventional crackers that are present in the superstructure. Depending on the nature of the energy carriers, the connections shown in the superstructure can represent either energy or mass flows.}\label{fig_superstructure}
\end{figure}

To address these challenges, in our previous work \cite{Naraghi2025}, we proposed a microgrid superstructure shown in Figure \ref{fig_superstructure} for steam cracking electrification. To account for the gradual and continuous electrification efforts in an ethylene plant, both conventional and electrified crackers are considered in the superstructure. Also, instead of solely relying on VRE to support 24/7 operation of electrified crackers, which requires massive battery storage that is impractically expensive \cite{AgrawalSiirola2023}, we adopt diverse energy sources as well as various energy storage and conversion systems in the superstructure. The microgrid consists of both dispatchable and non-dispatchable generators. The dispatchable generators are managed by the microgrid master controllers \cite{khodaei2013microgrid} obeying operational constraints such as generation capacity, ramping limits, and minimum on/off times, whereas non-dispatchable generators depend entirely on the availability of renewable energy sources. Together, we solved the day-ahead operational scheduling problem under different electrification levels and VRE profiles \cite{Naraghi2025}. We concluded that, given the current status of the power grid and renewable energy generation technologies, the process economics and sustainability of electrified steam cracking do not always favor higher decarbonization levels. To overcome this barrier, electricity from the main grid must be cleaner and cheaper, and energy storage costs per unit stored must decrease. Furthermore, both chemical and power systems stakeholders should seamlessly coordinate with each other to pursue joint optimization in operation. 

These findings motivate us to develop a joint optimization model that holistically integrates the operations of both electrified cracking microgrids and the bulk power system on a large scale. In this work, we propose a multi-agent unit commitment (MAUC) framework that simultaneously identifies the optimal operational schedule for electrified cracking microgrids and unit commitment and power dispatch schedules for the power grid connecting the microgrids. To efficiently solve the resulting large-scale MILP, we develop a two-stage solution technique based on Benders decomposition, which achieves a significant reduction in solution time without sacrificing accuracy. To illustrate the effectiveness of our approach, we investigate a case study on the electrification of steam crackers in the State of Texas. We survey all 26 ethylene plants in Texas and gathered their feedstocks, production capacities, and local wind and solar profiles. For each plant, we identify its associated bus on ACTIVSg2000, a 2000-bus synthetic grid benchmark test case for the Texas grid network \cite{Birchfield2017}. By bridging the gap between local operational strategies and system-level feasibility, this work marks the first systematic study to formulate and solve a joint optimization model for both electrified chemical plant microgrid and the bulk power grid. Results from this case study provide new perspectives on the interplay and systems-level dynamics between electrified chemical plants and the power system, thereby offering a systematic and holistic approach to evaluate the economic and operational benefits and barriers for industrial decarbonization via electrification.


\section{Literature Review}

\paragraph{Electrification of steam cracking.} Recent studies in electrification of steam cracking process heating primarily focus on new electrified reactor designs and plant-level decarbonization strategies \cite{Tijani2022,Mallapragada2023}. Earlier conceptual and comparative studies, such as Layritz et al. \cite{Layritz2021}, made clear that direct electrification can reduce process emissions only when the power supply is sufficiently low-carbon, a conclusion reinforced by subsequent optimization, life-cycle, and techno-economic studies. For example, Tiggeloven et al. \cite{Tiggeloven2023} examined electric ethylene production with grid power, batteries, and dedicated renewables and showed that operational flexibility can improve both cost and emissions performance, yet the abatement potential remains tightly linked to the carbon intensity of electricity. Likewise, Mynko et al. \cite{Mynko2023} found that electrified cracking concepts such as low-emission furnaces and rotodynamic reactors can outperform conventional steam cracking under renewable electricity, but may become counterproductive under fossil-dominated grids, while Cattry et al. \cite{Cattry2025} showed that electrified steam cracking can have higher life-cycle emissions than thermal steam cracking under current grid conditions and lower emissions only under cleaner future electricity scenarios. At the same time, much of the steam-cracking-specific work has focused on alternative reactor concepts and reactor design details, including modular electrical-resistance-heated reactors and internally heated wire-reactor configurations \cite{Balakotaiah2022,Agrawal2023,RodriguezGil2025}, plasma-based concepts \cite{Delikonstantis2019}, and rotodynamic cracking systems \cite{Mynko2023}. Although studies that move beyond reactor design toward flexible operation and renewable integration already exist \cite{Tiggeloven2023, Giannikopoulos2024}, they largely treat the electricity system through exogenous grid, storage, or renewable-supply assumptions rather than through endogenous coordination with a plant-level microgrid superstructure and the bulk power system. Overall, the existing literature suggests that electrified cracking reactor alone does not necessarily lead to lower emissions as it depends critically on where, when, and how electricity is generated and used. As a result, despite substantial progress in electrified reactor concepts and plant-level assessments, there remains a gap for a large-scale, integrated optimization framework that jointly considers the operational scheduling of electrified crackers, microgrids, and power systems.

\paragraph{Optimal operation of microgrids.} A microgrid is a localized electrical network that integrates at least one distributed energy resource and one load operates in both grid-connected and islanded modes \cite{TonSmith2012, khodaei2014resiliency}. Such dual capability enables microgrids to maintain continuity of supply during upstream disturbances in the main distribution system, thereby reducing the risk of load shedding and improving system resiliency \cite{khodaei2012microgrid, series2009microgrids}. In addition to reliability, microgrids offer a pathway for incorporating renewable energy resources into industrial operations, lowering greenhouse gas emissions and operational costs while enhancing power quality \cite{khodaei2012microgrid}. Early work on microgrids and distributed energy systems formulated day-ahead unit commitment and economic dispatch problems as deterministic MILP or MINLP models for systems with photovoltaics, gas turbines, boilers, combined heat and power, and batteries \cite{Yan2017,Kong2005,DiSomma2015,DiSomma2016}. Subsequent studies expanded this foundation to address renewable intermittency and operational uncertainty through scenario-based stochastic optimization, coordinated storage dispatch, and more recently robust and distributionally robust formulations for islanded and grid-connected microgrids \cite{NguyenCrow2016,Rahbar2015,Sandgani2017,Xie2022,Zhai2024,Jia2025,Abdelghany2024}. For industrial applications where thermal utilities are indispensable, joint dispatch models and decomposition-based coordination schemes have been developed for coupled electricity-steam or electricity-heating networks to improve computational tractability and privacy preservation in multi-stakeholder settings \cite{Li2016,Lin2017,Zheng2018,Lu2018,Chen2021}. In the chemical process systems engineering community, hydrogen- and ammonia-based microgrid studies have introduced increasingly process-aware formulations in which electrolyzers, Haber-Bosch synthesis, chemical storage, and dual-fuel power generation are co-optimized using combined design-and-scheduling models, hierarchical dynamic real-time optimization, rolling-horizon control, and data-driven distributionally robust optimization \cite{PALYS2020,Wang2021,Palys2019,Kong2024,Li2025}. Due to the uniqueness of process configuration, feedstock and product specifications, and operational challenges, each electrified chemical process requires a tailored microgrid superstructure design \cite{Bell2025}. For electrified cracking, however, no microgrid superstructure has been proposed until our recent work \cite{Naraghi2025, escape_25_microgrid}. Furthermore, prior work on microgrid operation mostly focuses on local unit commitment and resource scheduling without explicitly considering the broader interactions and coordination with the main power system, thereby posing an important research gap.

\paragraph{Unit commitment problem (UCP) of power systems.} The UCP is an $\mathcal{NP}$-hard problem that determines the optimal on/off scheduling of generators and their production levels subject to physical, operational, and economic constraints \cite{padhy2004unit}. Dynamic programming provided a systematic framework to handle various unit constraints, though it was hampered by exponential computational growth \cite{nieva2007cht, snyder2007dynamic}. Integer programming extended the branch-and-bound method to address practical scheduling problems \cite{dillon1978integer, lauer1982solution}, and MILP formulations later emerged as the dominant paradigm due to their ability to generate tighter approximations and near-optimal solutions efficiently \cite{carrion2006computationally, frangioni2008tighter, wu2011tighter}. Decomposition methods such as Lagrangian relaxation and Dantzig-Wolfe decomposition were also introduced to reduce computational burden. While these methods enabled tractable solutions for large problems, they often introduced suboptimality \cite{padhy2004unit, guan2018polynomial}. More recently, Benders decomposition has been proposed as a promising approach, dividing the UCP into a master problem and subproblems, thereby enabling parallel solution and reducing complexity while maintaining near-optimality. Despite these advances, the UCP remains a challenging problem to solve for large-scale power systems with significant renewable penetration \cite{Bendotti2019}.

\section{Multi-Agent Unit Commitment Problem Formulation}

We model the joint operational optimization of electrified steam cracking microgrids (assuming continuous, steady-state olefin production rates) and their connected power system as a multi-agent UCP, which is a common modeling choice for this type of problems \cite{nagata2002multi,yu2004solution,sharma2013multi}. Our goal is to minimize the total day-ahead operating costs for both the power system (denoted as $\mathrm{PS}$) and all chemical plants (denoted as $\mathrm{CP}$). We use an Independent System Operator (ISO) to manage power transfers, schedule generation, and ensure the optimal operation of both centralized and distributed energy resources. Rather than introducing a separate Distribution Market Operator (DMO) \cite{parhizi2016market}, the model assumes direct coordination between the ISO and microgrids. Furthermore, we utilize a hierarchical structure, where the main grid (ISO) is at the upper level and interacts with the microgrids at the lower level. Specifically, the ISO receives load demand from microgrids and determines unit schedules by solving the unit commitment problem. Meanwhile, the microgrid controller optimizes the lowest-cost schedule subject to local resources.

In the following formulations, we denote the set of all buses as $B$ and the (hourly) operational planning horizon as $T=\{1,2,\dots,24\}$. Let $K$ and $P$ be the complete set of generators and electrified steam cracking microgrids, respectively. Each bus $u\in B$ may be associated with a subset of generators and electrified cracking microgrids denoted as $K_u$ and $P_u$, respectively. In other words, $K = \cup_{u\in B} K_u$ and $P = \cup_{u\in B} P_u$. Furthermore, subset $B_u \subseteq B$ contains all buses connected to a given bus $u$. On the other hand, not all buses are connected to a microgrid. We define $B' \subseteq B$ as the subset of buses connected to at least one electrified cracking microgrid.


\subsection{Power System Formulation}\label{sec:ps}

The power system UCP formulation given by Problem \eqref{eq:psobj} minimizes the total operating costs of the power system, $z^{\mathrm{PS}}$:
\begin{equation}
    z^{\mathrm{PS}} = \sum_{t \in T} \sum_{k \in K} \left(c^{\mathrm{K}}_k x^{\mathrm{K}}_{k,t} + c^{\mathrm{d, K}}_k y^{\mathrm{K}}_{k,t} \right),\label{eq:psobj}
\end{equation}
where $y^{\mathrm{K}}_{k,t} \geq 0$ and $x^{\mathrm{K}}_{k,t} \in \{0,1\}$ denote the electricity dispatch and commitment decision variables for generator $k$ during hour $t$, respectively; and parameters $c^{\mathrm{d, K}}_k$ and $c^{\mathrm{K}}_k$ denote dispatch and commitment costs for generator $k$, respectively.

Constraint \eqref{eq:psminmaxcapacity} ensures that the power generation at each generator $k$ is bounded by its minimum capacity $\underline{P}^{\mathrm{K}}_k$ and maximum capacity $\overline{P}^{\mathrm{K}}_k$:
\begin{equation}
    \underline{P}^{\mathrm{K}}_k x^{\mathrm{K}}_{k,t} \leq  y^{\mathrm{K}}_{k,t} \leq \overline{P}^{\mathrm{K}}_k x^{\mathrm{K}}_{k,t}, \qquad \forall k \in K, t \in T\label{eq:psminmaxcapacity}
\end{equation}

Constraints within \eqref{eq:psupdown} where parameters $\mathrm{UT}^{\mathrm{K}}_k$ and $\mathrm{DT}^{\mathrm{K}}_k$ are the minimum up and down time, respectively; while $\pi^{\mathrm{K}}_{\mathrm{DT},k,t}$ and $\pi^{\mathrm{K}}_{\mathrm{UT},k,t}$ represent the up and down reserve deployment variables, respectively.
\begin{equation}\label{eq:psupdown}
    \begin{aligned}
        & -\pi^{\mathrm{K}}_{\mathrm{DT},k,t} \leq x^{\mathrm{K}}_{k,t} - x^{\mathrm{K}}_{k,t-1} \leq \pi^{\mathrm{K}}_{\mathrm{UT},k,t}, \qquad \forall k \in K,t \in T\backslash \{1\} \\
        & \sum_{\tau = t - \mathrm{UT}^{\mathrm{K}}_k + 1}^{t} \pi^{\mathrm{K}}_{\mathrm{UT},k,\tau} \leq x^{\mathrm{K}}_{k,t} \leq 1 - \sum_{\tau = t - \mathrm{DT}^{\mathrm{K}}_k + 1}^{t} \pi^{\mathrm{K}}_{\mathrm{DT},k,\tau} \qquad \forall k \in K,t \in T
    \end{aligned}
\end{equation}

Constraint \eqref{eq:psramping} ensures that each generator $k$ must obey its ramping up and down limits, denoted as $\mathrm{RU}^{K}_{k}$ and $\mathrm{RD}^{K}_{k}$, respectively:
\begin{equation}\label{eq:psramping}
    - \mathrm{RD}^{K}_{k} \leq y^{\mathrm{K}}_{k,t} - y^{\mathrm{K}}_{k,t-1} \leq \mathrm{RU}^{K}_{k} \qquad \forall k \in K, t \in T\backslash\{1\}
\end{equation}

Constraint \eqref{eq:phase} ensures that the power flow across link $u \in B$ and $\nu\in B_u$ is a function of the respective phase angles $\theta_{u,t}$ and $\theta_{\nu,t}$:
\begin{equation}
    \Gamma^{\mathrm{K}}_{u,\nu}\left(\theta_{u,t} - \theta_{\nu,t} \right) = f^{\mathrm{PS}}_{u,\nu,t}, \qquad \forall u \in B, \nu \in B_u, t \in T \label{eq:phase}
\end{equation}
where $f^{\mathrm{PS}}_{u,\nu,t}$ and $\Gamma^{\mathrm{K}}_{u,\nu}$ denote the hour $t$ power flow and phase angle-to-power conversion constant for line $u\nu$, respectively. In addition, Equation \eqref{eq:psphaseminmax} enforces network flow constraints:
\begin{equation}\label{eq:psphaseminmax}
    -F^{\mathrm{K}} \leq \Gamma^{\mathrm{K}}_{u,\nu}\left(\theta_{u,t} - \theta_{\nu,t} \right) \leq F^{\mathrm{K}}, \qquad \forall u \in B, \nu \in B_u, t \in T 
\end{equation}
where $F^{\mathrm{K}}$ denotes the transmission capacity limit.

Finally, Constraint \eqref{eq:psflowbalance} is the bus-wise power balance indicating that, for each bus $u\in B$, the total power demand of all connected chemical plant microgrids $\zeta_{u,t}$ and other users $\delta_{u,t}$ must be balanced by the sum of power generated by the attached generators and the net power flow into $u$:
\begin{equation}
    \sum_{k \in K_u}y^{\mathrm{K}}_{k,t} + \eta^{\mathrm{K}}_{u,t} - \delta_{u,t} - \zeta_{u,t} = \sum_{\nu \in B_u}  \Gamma^{\mathrm{K}}_{u,\nu}\left(\theta_{u,t} - \theta_{\nu,t} \right), \qquad \forall u \in B, t \in T \label{eq:psflowbalance}
\end{equation}
where $\eta^{\mathrm{K}}_{u,t}$ is the renewable power generation at bus $u$. Note that the RHS of Constraint \eqref{eq:psflowbalance} gives the difference between the total power flow leaving bus $u$ and the total power flow coming to $u$.

\subsection{Electrified Steam Cracking Microgrid Planning Formulation}\label{sec:cp}

The electrified steam cracking microgrid formulation given by Problem \eqref{eq:cpobj} minimizes the total operating costs of all microgrids, $z^{\mathrm{CP}}$, which account for fuel (natural gas) expenses, electricity costs from the grid, energy generation costs, startup and shutdown costs, and hydrogen storage costs:
\begin{equation}\label{eq:cpobj}
    \begin{aligned}
        z^{\mathrm{CP}}= \sum_{p \in P}\sum_{t \in T} & \Bigg(c^{\mathrm{F}}_p f^{\mathrm{F}}_{\mathrm{CC},p,t} + \sum_{i \in I_p} c^{\mathrm{F}}_p f^{\mathrm{F,I}}_{i,p,t} + c^{\mathrm{G}}_{p,t}p^{\mathrm{G}}_{p,t} + c^{\mathrm{FC}}_p p^{\mathrm{FC}}_{p,t} + c^{\mathrm{EL}}_p f^{\ce{H2}}_{\mathrm{EL},p,t} + c^{\mathrm{HS}}_p f^{\ce{H2}}_{\mathrm{HS},p,t} \\
        &+ \sum_{i \in I_p} \Big[c^{\mathrm{I}}_{i,p}p^{\mathrm{I}}_{i,p,t} + c^{\mathrm{I, SU}}_{i,p} \max \{0, x^{\mathrm{I}}_{i,p,t} - x^{\mathrm{I}}_{i,p,t-1}\} + c^{\mathrm{I, SD}}_{i,p} \max \{0, x^{\mathrm{I}}_{i,p,t-1} - x^{\mathrm{I}}_{i,p,t} \} \Big] \Bigg),
    \end{aligned}
\end{equation}
where parameters $c^{\mathrm{F}}_p$, $c^{\mathrm{G}}_{p,t}$, $c^{\mathrm{FC}}_p$, $c^{\mathrm{EL}}_p$, $c^{\mathrm{HS}}_p$, $c^{\mathrm{I}}_{i,p}$, $c^{\mathrm{I, SU}}_{i,p}$, and $c^{\mathrm{I, SD}}_{i,p}$ represent the plant $p$'s fuel cost, electricity market price at time $t$, fuel cell operating cost, PEM electrolyzer operating cost, hydrogen storage unit cost, local fuel-based generator operating cost, local fuel-based generator startup ($\mathrm{SU}$) and shutdown ($\mathrm{SD}$) costs, respectively. Non-negative continuous variables $f^{\mathrm{F}}_{\mathrm{CC},p,t}$, $f^{\mathrm{F, I}}_{i,p,t}$, $p^{\mathrm{G}}_{p,t}$, $p^{\mathrm{FC}}_{p,t}$, $f^{\ce{H2}}_{\mathrm{EL},p,t}$, $f^{\ce{H2}}_{\mathrm{HS},p,t}$, $p^{\mathrm{I}}_{i,p,t}$ denote plant $p$'s fuel consumption rate by conventional fuel-based crackers (abbreviated as $\mathrm{CC}$), fuel consumption rate by local fuel-based generators, grid power usage by electrified crackers (abbreviated as $\mathrm{EC}$), power generation rate by fuel cell units, \ce{H2} production rate from PEM electrolyzers, \ce{H2} storage level, and power generation rate from local fuel-based generators, respectively. Binary variable $x^{\mathrm{I}}_{i,p,t}$ represents the commitment status during hour $t$ of local fuel-based generator $i$ belonging to plant $p$ ($x^{\mathrm{I}}_{i,p,0}$ indicates initial commitment status). Set $I_p$ is a set of local fuel-based generators associated with microgrid $p$. Note that the maximization operations in calculating the local generator startup and shutdown costs in Problem \eqref{eq:cpobj} can be easily linearized.

Constraint \eqref{eq:cpenergy} describes the different sources of energy supplied to conventional and electrified crackers. 
\begin{equation}\label{eq:cpenergy}
    \begin{aligned}
        & \dot{\mathrm{Q}}^{\ce{CH4}} \left( f^{\mathrm{F}}_{\mathrm{CC},p,t} +  f^{\ce{CH4}}_{\mathrm{Sep},\mathrm{CC},p,t}\right) + \dot{\mathrm{Q}}^{\ce{H2}} \left( f^{\ce{H2}}_{\mathrm{Sep,CC}, p,t} +  f^{\ce{H2}}_{\mathrm{HS,CC},p,t} + f^{\ce{H2}}_{\mathrm{EL,CC},p,t} \right) = P^{\mathrm{CC}}_p, \\
        & \sum_{i \in I_p} p^{\mathrm{I}}_{\mathrm{EC},i,p,t} + p^\mathrm{RE}_{\mathrm{EC},p,t} + p^{\mathrm{G}}_{\mathrm{EC},p,t} + p^\mathrm{ESS}_{\mathrm{EC},p,t} + p^\mathrm{FC}_{\mathrm{EC},p,t} = P^{\mathrm{EC}}_p,
    \end{aligned}
    \forall p \in P, t \in T
\end{equation}
in which, for conventional crackers in microgrid $p$, parameter $\dot{\mathrm{Q}}^{\ce{CH4}}$ (resp. $\dot{\mathrm{Q}}^{\ce{H2}}$) is the lower heating value of natural gas/methane fuel (resp. \ce{H2} fuel), variable $f^{\mathrm{F}}_{\mathrm{CC},p,t}$ is the natural gas fuel consumption from conventional crackers, variable $f^{\ce{CH4}}_{\mathrm{Sep},\mathrm{CC},p,t}$ (resp. $f^{\ce{H2}}_{\mathrm{Sep},\mathrm{CC},p,t}$) is the flow rate of \ce{CH4} (resp. \ce{H2}) recovered from the separation units and sent back to $\mathrm{CC}$, $f^{\ce{H2}}_{\mathrm{HS},\mathrm{CC},p,t}$ ($f^{\ce{H2}}_{\mathrm{EL},\mathrm{CC},p,t}$) stands for the flow rate of \ce{H2} sent from \ce{H2} storage unit (resp. electrolyzer) to $\mathrm{CC}$, and parameter $P^{\mathrm{CC}}_p$ is the total heat power required for all conventional crackers (which can be determined by cracking reaction kinetics, plant size, and feedstock type; see our earlier work \cite{Naraghi2025} for details). Meanwhile, for electrified crackers in microgrid $p$, the total electrical power required $P^{\mathrm{EC}}_p$ (which can also be determined a priori following \cite{Naraghi2025}) is met by the electrical power received from 1) local fuel-based generators $p^{\mathrm{I}}_{\mathrm{EC},i,p,t}$, 2) non-dispatchable units (solar and wind) $p^\mathrm{RE}_{\mathrm{EC},p,t}$, 3) grid power from the main power system $p^\mathrm{G}_{\mathrm{EC},p,t}$, 4) local energy storage system (battery units) $p^\mathrm{ESS}_{\mathrm{EC},p,t}$, and 5) local hydrogen fuel cell units $p^\mathrm{FC}_{\mathrm{EC},p,t}$.

Constraint \eqref{eq:cph2balance} states that 1) the light cracking product (i.e., \ce{H2} and \ce{CH4}) from $\mathrm{CC}$ (resp. $\mathrm{EC}$) sent to downstream separation units for recovery, namely $f^{\ce{H2}/\ce{CH4}}_{\mathrm{CC,Sep},p,t}$ (resp. $f^{\ce{H2}/\ce{CH4}}_{\mathrm{EC,Sep},p,t}$), cannot exceed its generation rate (which is a parameter determined a priori from \cite{Naraghi2025}), $F^{\mathrm{CC,Mix}}_p$ (resp. $F^{\mathrm{EC,Mix}}_p$), and 2) \ce{H2} leaving the separation unit can either be sent back to the $\mathrm{CC}$ as fuel, or can be sent to fuel cells and \ce{H2} storage unit, while \ce{CH4} leaving the separation unit will be sent back to the the $\mathrm{CC}$ as fuel.
\begin{equation}\label{eq:cph2balance}
    \begin{aligned}
        & 0 \leq f^{\ce{H2}/\ce{CH4}}_{\mathrm{CC,Sep},p,t}\leq F^{\mathrm{CC,Mix}}_p,\\
        & 0 \leq f^{\ce{H2}/\ce{CH4}}_{\mathrm{EC,Sep},p,t}\leq F^{\mathrm{EC,Mix}}_p, \\
        & f^{\ce{CH4}}_{\mathrm{Sep,CC}, p,t} = R^{\ce{CH4}}r_p^{\ce{CH4}}\left(f^{\ce{H2}/\ce{CH4}}_{\mathrm{CC,Sep},p,t} + f^{\ce{H2}/\ce{CH4}}_{\mathrm{EC,Sep},p,t} \right), \\
        & f^{\ce{H2}}_{\mathrm{Sep,CC}, p,t} + f^{\ce{H2}}_{\mathrm{Sep,FC},p,t} + f^{\ce{H2}}_{\mathrm{Sep,HS}, p,t} = R^{\ce{H2}}r_p^{\ce{H2}}\left(f^{\ce{H2}/\ce{CH4}}_{\mathrm{CC,Sep},p,t} + f^{\ce{H2}/\ce{CH4}}_{\mathrm{EC,Sep},p,t} \right),
    \end{aligned}
    \qquad \forall p \in P, t \in T
\end{equation}
where variables $f^{\ce{CH4}}_{\mathrm{Sep,CC}, p,t}$, $f^{\ce{H2}}_{\mathrm{Sep,CC}, p,t}$, $f^{\ce{H2}}_{\mathrm{Sep,FC}, p,t}$, and $f^{\ce{H2}}_{\mathrm{Sep,HS}, p,t}$ denote \ce{CH4} product sent from separation unit to $\mathrm{CC}$, \ce{CH4} product sent from separation unit to $\mathrm{CC}$, fuel cells, and hydrogen storage unit, respectively. Parameter $r_p^{\ce{H2}}$ (resp. $r_p^{\ce{CH4}}$) is the mass flow ratio of \ce{H2} (resp. \ce{CH4}) in the light cracking product (i.e., $r_p^{\ce{H2}} + r_p^{\ce{CH4}} = 1$); and parameter $R^{\ce{H2}}$ (resp. $R^{\ce{CH4}}$) is the \ce{H2} (resp. \ce{CH4}) recovery in downstream separation unit.

Constraint \eqref{eq:cph2hs} models \ce{H2} balance around the storage unit in microgrid $p$: any \ce{H2} deficit or excess will result in a decrease or increase in \ce{H2} storage level, which is bounded by the storage capacity $\mathrm{HSC}_p$. 
\begin{equation}\label{eq:cph2hs}
    \begin{aligned}
        & f^{\ce{H2}}_{\mathrm{HS},p,t-1} + f^{\ce{H2}}_{\mathrm{Sep,HS},p,t} + f^{\ce{H2}}_{\mathrm{EL,HS},p,t} - f^{\ce{H2}}_{\mathrm{HS,CC},p,t} - f^{\ce{H2}}_{\mathrm{HS,FC},p,t} =  f^{\ce{H2}}_{\mathrm{HS},p,t}, \\
        & 0 \leq f^{\ce{H2}}_{\mathrm{HS},p,t} \leq \mathrm{HSC}_p,\\
    \end{aligned}
    \qquad \forall p \in P, t \in T
\end{equation}
where $f^{\ce{H2}}_{\mathrm{EL,HS},p,t}$ and $f^{\ce{H2}}_{\mathrm{HS,FC},p,t}$ stand for the flow rates of \ce{H2} from electrolyzer to storage and from \ce{H2} storage to the fuel cell unit, respectively; and $f^{\ce{H2}}_{\mathrm{HS},p,0}$ denotes the initial hydrogen storage level at the beginning of the planning horizon ($t=0$).

Constraint \eqref{eq:cph2el} models \ce{H2} production from the PEM electrolyzers in microgrid $p$, $f^{\ce{H2}}_{\mathrm{EL},p,t}$, bounded by the production capacity $\mathrm{PEL}_p$.
\begin{equation} \label{eq:cph2el}
    \begin{aligned}
        & f^{\ce{H2}}_{\mathrm{EL},p,t} = f^{\ce{H2}}_{\mathrm{EL,CC},p,t} + f^{\ce{H2}}_{\mathrm{EL,FC},p,t} + f^{\ce{H2}}_{\mathrm{EL,HS},p,t}, \\
        & f^{\ce{H2}}_{\mathrm{EL},p,t} = \frac{\eta^{\mathrm{EL}}}{P^{\ce{H2}}} \left(\sum_{i \in I_p} p^{\mathrm{I}}_{\mathrm{EL},i,p,t} + p^{\mathrm{RE}}_{\mathrm{EL},p,t} + p^{\mathrm{G}}_{\mathrm{EL},p,t} + p^{\mathrm{ESS}}_{\mathrm{EL},p,t} \right), \\
        & 0 \leq f^{\ce{H2}}_{\mathrm{EL},p,t} \leq \mathrm{PEL}_p,
    \end{aligned}
    \qquad \forall p \in P, t \in T
\end{equation}
where variable $f^{\ce{H2}}_{\mathrm{EL,FC},p,t}$ is the \ce{H2} flow rate from electrolyzer to fuel cell unit, variable $p^{\mathrm{G}}_{\mathrm{EL},p,t}$ is electrolyzer power input from the grid, parameter $\eta^{\mathrm{EL}}$ denotes electrical efficiency of PEM electrolyzer (e.g., 73.6\% \cite{pemefficiency}), and parameter $P^{\ce{H2}}$ corresponds to the electrolysis energy requirement with no efficiency loss (39.4 MWhr/ton of \ce{H2}). 

Now proceeding to fuel cell units, Constraint \eqref{eq:cpfs} relates fuel cell commitment status (indicated by binary variable $x^{\mathrm{FC}}_{p,t}$) with its power generation and distribution.
\begin{equation} \label{eq:cpfs}
\begin{aligned}
    & \underline{P}^{\mathrm{FC}}_p x^{\mathrm{FC}}_{p,t} \leq p^{\mathrm{FC}}_{p,t} \leq \overline{P}^{\mathrm{FC}}_p x^{\mathrm{FC}}_{p,t},\\
    & p^{\mathrm{FC}}_{p,t} = \eta^{\mathrm{FC}} \dot{\mathrm{Q}}^{\ce{H2}} \left( f^{\ce{H2}}_{\mathrm{HS,FC},p,t} + f^{\ce{H2}}_{\mathrm{EL,FC},p,t} + f^{\ce{H2}}_{\mathrm{Sep,FC},p,t} \right),\\
    & p^{\mathrm{FC}}_{p,t} = p^{\mathrm{FC}}_{\mathrm{EC},p,t}+ p^{\mathrm{FC}}_{\mathrm{ESS},p,t},
\end{aligned}
    \qquad \forall p \in P, t \in T
\end{equation}
where variables $p^{\mathrm{FC}}_{\mathrm{EC},p,t}$ and $p^{\mathrm{FC}}_{\mathrm{ESS},p,t}$ denote power flow from fuel cell units to electrified crackers and energy storage system (batteries), respectively; and parameters $\eta^{\mathrm{FC}}$, $\underline{P}^{\mathrm{FC}}_p$, and $\overline{P}^{\mathrm{FC}}_p$ correspond to energy efficiency of fuel cell (e.g., 65\% \cite{fuelcellefficiency}), as well as the lower and upper generation capacities of the fuel cell, respectively.

Constraint \eqref{eq:cpgrid} models the distribution of grid across onsite local units.
\begin{equation}\label{eq:cpgrid}
    \begin{aligned}
        & 0 \leq p^{\mathrm{G}}_{p,t} \leq \overline{P}^{\mathrm{G}}_{p,t},\\
        & p^{\mathrm{G}}_{p,t} = p^{\mathrm{G}}_{\mathrm{EC},p,t} + p^{\mathrm{G}}_{\mathrm{ESS},p,t} + p^{\mathrm{G}}_{\mathrm{EL},p,t},
    \end{aligned}
    \qquad \forall p \in P, t \in T
\end{equation}
where variable $p^{\mathrm{G}}_{\mathrm{ESS},p,t}$ is the grid power directly stored in batteries, and $\overline{P}^{\mathrm{G}}_{p,t}$ is the maximum capacity grid power can be transferred to plant $p$.

Constraint \eqref{eq:cprenewable} models the distribution of local renewable power generation to electrified crackers, battery storage, and electrolyzer units. 
\begin{equation}\label{eq:cprenewable}
    p^{\mathrm{RE}}_{\mathrm{EC},p,t} + p^{\mathrm{RE}}_{\mathrm{ESS},p,t} + p^{\mathrm{RE}}_{\mathrm{EL},p,t} = P^{\mathrm{PV}}_{p,t} + P^{\mathrm{WT}}_{p,t}, \qquad \forall p \in P, t \in T
\end{equation}
where parameters $P^{\mathrm{PV}}_{p,t}$ and $P^{\mathrm{WT}}_{p,t}$ correspond to hourly power generation from local PV solar panels and wind turbines installed in plant $p$.

The operational constraints related to local fuel-based (dispatchable) generators are formulated in a similar manner as the power system UCP discussed in Section \ref{sec:ps}. They consist of power capacity limits, fuel consumption, power distribution, ramp-up and ramp-down limits, and minimum up/down time requirements. 
\begin{equation}
    \begin{aligned}
        &\underline{P^{\mathrm{I}}}_{i,p} x^{\mathrm{I}}_{i,p,t} \leq p^{\mathrm{I}}_{i,p,t} \leq \overline{P^{\mathrm{I}}}_{i,p} x^{\mathrm{I}}_{i,p,t},\\
        & p^{\mathrm{I}}_{i,p,t} = \eta^{\mathrm{I}}_{i} \dot{\mathrm{Q}}^{\ce{CH4}} f^{\mathrm{F,I}}_{i,p,t}, \\
        & p^{\mathrm{I}}_{i,p,t} = p^{\mathrm{I}}_{\mathrm{EC},i,p,t} + p^{\mathrm{I}}_{\mathrm{ESS},i,p,t} + p^{\mathrm{I}}_{\mathrm{EL},i,p,t},\\
        & -\mathrm{RD}^{\mathrm{I}}_{i,p}\leq p^{\mathrm{I}}_{i,p,t} - p^{\mathrm{I}}_{i,p,t-1} \leq \mathrm{RU}^{\mathrm{I}}_{i,p},\\
        & x^{\mathrm{I}}_{i,p,t} - x^{\mathrm{I}}_{i,p,t-1} \leq x^{\mathrm{I}}_{i,p,\tau}, \quad t+1 \leq \tau \leq \min \{|T|,t+\mathrm{UT}_{i,p}-1\}\\
        & x^{\mathrm{I}}_{i,p,t-1} - x^{\mathrm{I}}_{i,p,t} \leq 1- x^{\mathrm{I}}_{i,p,\tau}, \quad t+1 \leq \tau \leq \min \{|T|,t+\mathrm{DT}_{i,p}-1\}
    \end{aligned}
    \quad \forall i \in I_p, p \in P, t \in T
\end{equation}
where parameter $\eta^{\mathrm{I}}_{i}$ denotes the efficiency of gas turbine combined cycle (e.g., 60\% \cite{gtccefficiency}), $\mathrm{RU}^{\mathrm{I}}_{i,p}$ (resp. $\mathrm{RD}^{\mathrm{I}}_{i,p}$) is ramping up (resp. ramping down) limit for local generator $i$, and $\mathrm{UT}_{i,p}$ (resp. $\mathrm{DT}_{i,p}$) is the minimum up (resp. down) time for generator $i$.

Constraint \eqref{eq:cpess} models the power balance and charging/discharging mechanisms of the energy storage system made of batteries. For simplicity, we assume that every battery cell behaves the same. It accounts for the minimum charging time $\mathrm{MC}^{\mathrm{ESS}}_{p}$ and discharging time $\mathrm{MD}^{\mathrm{ESS}}_{p}$, battery storage capacity $\mathrm{ESC}_{p}$ and its minimum (resp. maximum) charging and discharging limits, $\underline{P}^{\mathrm{ESS,Ch}}_{p,t}$ and $\underline{P}^{\mathrm{ESS,Dch}}_{p,t}$ (resp. $\overline{P}^{\mathrm{ESS,Ch}}_{p,t}$ and $\overline{P}^{\mathrm{ESS,Dch}}_{p,t}$).
\begin{equation}\label{eq:cpess}
    \begin{aligned}
        & x^{\mathrm{ESS,Ch}}_{p,t} + x^{\mathrm{ESS,Dch}}_{p,t} \leq 1,\\
        & \underline{P}^{\mathrm{ESS,Ch}}_{p,t} x^{\mathrm{ESS,Ch}}_{p,t} \leq p^{\mathrm{Ch}}_{p,t} \leq \overline{P}^{\mathrm{ESS,Ch}}_{p,t} x^{\mathrm{ESS,Ch}}_{p,t} \\
        & \underline{P}^{\mathrm{ESS,Dch}}_{p,t} x^{\mathrm{ESS,Dch}}_{p,t} \leq p^{\mathrm{Dch}}_{p,t} \leq \overline{P}^{\mathrm{ESS,Dch}}_{p,t} x^{\mathrm{ESS,Dch}}_{p,t}, \\
        & x^{\mathrm{ESS,Ch}}_{p,t} - x^{\mathrm{ESS,Ch}}_{p,t-1} \leq x^{\mathrm{ESS,Ch}}_{p,\tau},\; t+1 \leq \tau \leq \min \{|T|,t+\mathrm{MC}^{\mathrm{ESS}}_{p}-1\}, \\
        & x^{\mathrm{ESS,Dch}}_{p,t} - x^{\mathrm{ESS,Dch}}_{p,t-1} \leq x^{\mathrm{ESS,Dch}}_{p,\tau}, \; t+1 \leq \tau \leq \min \{|T|,t+\mathrm{MD}^{\mathrm{ESS}}_{p}-1\}, \\
        & p^{\mathrm{ESS}}_{p,t} = p^{\mathrm{ESS}}_{p,t-1} + p^{\mathrm{Ch}}_{p,t} - p^{\mathrm{Dch}}_{p,t}, \\
        & 0 \leq p^{\mathrm{ESS}}_{p,t} \leq \mathrm{ESC}_{p}, \\
        & p^{\mathrm{Ch}}_{p,t} = p^{\mathrm{G}}_{\mathrm{ESS},p,t} + \sum_{i \in I_p}p^{\mathrm{I}}_{\mathrm{ESS},i,p,t} + p^{\mathrm{RE}}_{\mathrm{ESS},p,t} + p^{\mathrm{FC}}_{\mathrm{ESS},p,t}, \\
        & p^{\mathrm{Dch}}_{p,t} = p^{\mathrm{ESS}}_{\mathrm{EC},p,t} + p^{\mathrm{ESS}}_{\mathrm{EL},p,t}, \\
    \end{aligned}
    \qquad \forall p\in P, t \in T
\end{equation}
in which binary variables $x^{\mathrm{ESS,Ch}}_{p,t}$ and $x^{\mathrm{ESS,Dch}}_{p,t}$ specify the battery charging and discharging status, respectively; and continuous variables $p^{\mathrm{ESS}}_{p,t}$, $p^{\mathrm{Ch}}_{p,t}$, and $p^{\mathrm{Dch}}_{p,t}$ respectively represent the battery stored energy level, battery charging rate, and discharging rate for plant $p$ during time $t$.

\subsection{Coupling Constraint and the Joint Optimization Problem}

The power system UCP and the microgrid operational planning problem are linked together via a coupling constraint. Constraint \eqref{eq:coupling} ensures that the total power transferred through a bus to its connected electrified cracking microgrids must exactly match the sum of their grid power demands.
\begin{equation}\label{eq:coupling}
    \zeta_{u,t} = \sum_{p \in P_u}p^{\mathrm{G}}_{p,t}, \qquad \forall u \in B', t \in T
\end{equation}

Overall, Equations \eqref{eq:psobj} through \eqref{eq:coupling} form the original multi-agent UCP, denoted as (P), whose objective is to minimize the total operating costs $z^{\mathrm{PS}} + z^{\mathrm{CP}}$.

\section{Solution Approach}

\begin{figure}[ht!]
    \centering
    \includegraphics[width = 0.75\textwidth]{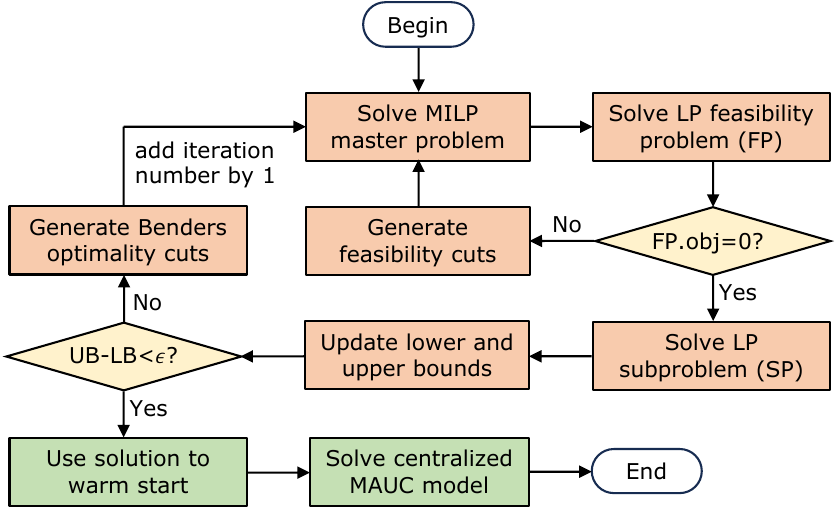}
    \caption{Two-stage procedure to solve the large-scale MILP multi-agent UCP.}
    \label{algorithm}
\end{figure}

Given the scale of the multi-agent UCP, solving (P) directly as a large-scale MILP would be computationally expensive. To improve computational efficiency, we propose a two-stage solution methodology based on Benders' decomposition and MILP refinement techniques illustrated in Figure \ref{algorithm}. Leveraging the hierarchical structure of the problem, where the power system is at the upper level and interacts with the electrified cracking microgrids located at the lower level, we decompose the original problem (P) into a master problem featuring the power system and multiple subproblems, each associated with a bus connected to chemical plants. To speed up the solution process, in the first stage, we solve the relaxed subproblems by removing their integrality constraints. The decomposition process iterates by generating feasibility and optimality cuts from the subproblems and incorporating them into the master problem, progressively refining the solution until convergence is achieved. Next, the near-optimal solution obtained from Benders' decomposition will serve as a warm start for the second stage, where the full, original MILP problem (P) is solved to ensure feasibility and optimality.

\subsection{First Stage: Benders' Decomposition}

\paragraph{MILP master problem (MP).} We formulate the master problem as follows:
\begin{equation}
    \begin{aligned}
        (\mathrm{MP}): \quad \min \; & z^{\mathrm{PS}} + \sum_{u \in B'}\phi_u \\
        & \text{Constraints \eqref{eq:psminmaxcapacity} through \eqref{eq:psflowbalance}}, \\
        & \phi^{\mathrm{LB}}_u \leq \phi_u, \quad \forall u \in B'
    \end{aligned}
\end{equation}
where \(\phi_u\) is an underestimator for the objective function value of subproblem ($\mathrm{SP}_u$), and \(\phi_u^{\mathrm{LB}}\) is a proper lower bound for \(\phi_u\) to prevent unboundedness.

\paragraph{Feasibility problem and feasibility cut.} After solving (MP), we obtain the optimal decision variables, including the power withdrawn from the power system to all connected chemical plant microgrids $\zeta^{*}_{u,t}$ for every $u \in B'$. However, not all $\zeta^{*}_{u,t}$ would satisfy the coupling constraint of \eqref{eq:coupling}. To identify the associated buses, we solve the following feasibility problem (FP) for each bus $u$, which is a linear program as all binary variables are relaxed.
\begin{equation}
    \begin{aligned}
        (\mathrm{FP}_u): \quad \min & \sum_{t \in T} \left(v^{+}_{u,t} + v^{-}_{u,t} \right) \\
         & \text{Constraints \eqref{eq:cpenergy} through \eqref{eq:cpess} (binary variables are relaxed)}, \\
         & \sum_{p \in P_u} p^{\mathrm{G}}_{p,t} - \zeta^{*}_{u,t} + v^{+}_{u,t} - v^{-}_{u,t} = 0. \quad \forall t \in T
    \end{aligned}
\end{equation}

If the optimal objective function value of ($\mathrm{FP}_u$) is 0, then $\zeta^{*}_{u,t}$ is feasible for bus $u$. Otherwise, a feasibility cut is added to (MP) as a constraint:
\begin{equation}\label{eq:feasibilitycut}
    \hat{\lambda}_{u,t}\left(\zeta_{u,t} - \zeta^{*}_{u,t} \right) + (v^{+}_{u,t}+ v^{-}_{u,t})^{*} \leq 0, \qquad \forall u \in B', t \in T
\end{equation}
where $\hat{\lambda}_{u,t}$ is the optimal dual value and $(v^{+}_{u,t}+ v^{-}_{u,t})^{*}$ is the optimal objective function value of ($\mathrm{FP}_u$).

\paragraph{Relaxed subproblems and Benders' optimality cut.} After solving (MP) and (FP) iteratively and obtaining a feasible \(\zeta^{*}\) for all \(u \in B'\), we will solve the relaxed subproblem (\(\mathrm{SP}_{u}\)) for each \(u \in B'\).
\begin{equation}
\begin{aligned}
    (\mathrm{SP}_{u}): \quad \min \; & z^{\mathrm{CP}}_{u}\\
    & \text{Constraints \eqref{eq:cpenergy} through \eqref{eq:cpess} (binary variables are relaxed)}, \\
    & \sum_{p \in P_u}p^{\mathrm{G}}_{p,t} - \zeta^{*}_{u,t} = 0, \qquad \forall t \in T
\end{aligned}
\end{equation}
in which the objective function $z^{\mathrm{CP}}_{u}$ is modified from $z^{\mathrm{CP}}$ in \eqref{eq:cpobj} by summing over $p \in P_u$.

If none of the ($\mathrm{SP}_{u}$) is unbounded, then Benders' optimality cuts of \eqref{eq:benderscut} are generated and added to the (MP).
\begin{equation}\label{eq:benderscut}
    {z^{\mathrm{CP}}_{u}}^* + \sum_{t \in T}\lambda_{u,t}\left(\zeta_{u,t}-\zeta^{*}_{u,t}\right) \leq \phi_u, \qquad \forall u \in B' 
\end{equation}
where $\lambda_{u,t}$ is the optimal dual value and ${z^{\mathrm{CP}}_{u}}^*$ is the optimal objective function value of ($\mathrm{SP}_u$).

\paragraph{Convergence.} After every iteration, the lower bound \(z^{\mathrm{LB}} = z^{\mathrm{PS}} + \sum_{u \in B'}\phi_u\) and the upper bound \(z^{\mathrm{UB}} = z^{\mathrm{PS}} + \sum_{u \in B'}z^{\mathrm{CP}}_u\) are updated and checked for convergence based on the relative gap \(\epsilon\) (usually \(\leq 10^{-4}\)). If 
\begin{equation}
    \frac{|z^{\mathrm{UB}} - z^{\mathrm{LB}}|}{|z^{\mathrm{UB}}|} \leq \epsilon,
\end{equation}
then Benders' decomposition converges, and the first-stage decomposition algorithm is terminated.

\subsection{Second Stage: Solving the Original Multi-agent UCP}

The optimal solution obtained from the first stage may not be optimal for the original multi-agent UCP, as the integrality constraints in the electrified cracking microgrid formulation are relaxed. However, it serves as a high-quality initial guess (near optimal), as the integrality constraints for the power system model are preserved throughout Benders' decomposition, ensuring a well-structured solution. Thus, we use the first-stage output to warm start the full MILP using commercial solvers such as Gurobi or CPLEX. Specifically, if the optimal solution from the first stage is feasible for the original MILP, we initialize both binary and continuous variables in both power system and microgrid formulations. Otherwise, we only initialize the binary variables associated with the power system unit commitment problem. We show that this significantly reduces the overall solution time compared to solving the full problem from scratch.

\section{A Full-scale Case Study: Electrification of Ethylene Production in Texas}

\subsection{Problem Setup}

To investigate the operational performance of our proposed electrification strategy and centralized optimization methodology, we conduct the first full-scale realistic case study for the joint optimization of the Texas grid and 26 ethylene plants in Texas. As shown in Table \ref{tab:plants}, we meticulously surveyed these ethylene plants and gathered their locations, feedstock compositions, and ethylene production capacities from various news sources. Based on this information, the minimum energy requirement of each ethylene plant can be determined using a dynamic optimization model presented in our earlier work \cite{Naraghi2025}, which incorporates steam cracking reaction kinetics and reactor design features. We remark that, by the time this study was performed (Spring 2025), Table \ref{tab:plants} represents the most accurate and complete information that includes all the ethylene plants in operation in Texas. It is worth noting that, although some of the values presented in Table \ref{tab:plants} may change since they were curated, the scale of the joint optimization problem, the methodologies used to analyze the problem, as well as the key findings and major trends observed in the results should not change.

\begin{table}[ht!]
\centering
\resizebox{\textwidth}{!}{
\begin{tabular}{clllcccccc}
\toprule
\bf Plant & \bf Company & \bf Site name & \bf County & \bf Ethylene prod. & \bf Ethane & \bf Propane & \bf Naphtha & \bf Gas oil & \bf Bus \\
\bf \# & & & & \bf (Mton/yr) & \bf (\%) & \bf (\%) & \bf (\%) & \bf (\%) & \bf ID \\
\midrule
1  & Chevron Phillips Chemical & Sweeny           & Brazoria     & 2.25   & 65  & 35   & 0   & 0   & 1612 \\
2  & Dow Chemical Co.          & Freeport         & Brazoria     & 6      & 50  & 50   & 0   & 0   & 1837 \\
3  & INEOS Olefins and Polymers& Chocolate Bayou  & Brazoria     & 1.89   & 100 & 0    & 0   & 0   & 1494 \\
4  & Formosa Plastics Corp. USA& Point Comfort    & Calhoun      & 1.2    & 100 & 0    & 0   & 0   & 1467 \\
5  & Formosa Plastics Corp. USA& Point Comfort    & Calhoun      & 3.5    & 45  & 25   & 15  & 15  & 1467 \\
6  & Eastman Chemical Co.      & Longview         & Gregg        & 0.64   & 30  & 55   & 15  & 0   & 1987 \\
7  & Chevron Phillips Chemical & Cedar Bayou      & Harris       & 1.5    & 100 & 0    & 0   & 0   & 1687 \\
8  & Chevron Phillips Chemical & Cedar Bayou      & Harris       & 2.2    & 25  & 45   & 30  & 0   & 1687 \\
9  & ExxonMobil Chemical Co.   & Baytown          & Harris       & 1.5    & 100 & 0    & 0   & 0   & 1691 \\
10 & ExxonMobil Chemical Co.   & Baytown          & Harris       & 2.3    & 0   & 0    & 100 & 0   & 1691 \\
11 & LyondellBasell Industries & Channelview      & Harris       & 1.93   & 5   & 10   & 40  & 45  & 1775 \\
12 & LyondellBasell Industries & La Porte         & Harris       & 1.15   & 80  & 20   & 0   & 0   & 1539 \\
13 & Shell Chemicals Ltd.      & Deer Park        & Harris       & 0.961  & 10  & 0    & 60  & 30  & 1513 \\
14 & BASF/Total                & Port Arthur      & Jefferson    & 1      & 40  & 40   & 20  & 0   & 1648 \\
15 & Bayport Polymers, LLC     & Port Arthur      & Jefferson    & 1      & 100 & 0    & 0   & 0   & 1672 \\
16 & Chevron Phillips Chemical & Port Arthur      & Jefferson    & 0.855  & 0   & 0    & 100 & 0   & 1648 \\
17 & ExxonMobil Chemical Co.   & Beaumont         & Jefferson    & 0.826  & 70  & 24   & 6   & 0   & 1648 \\
18 & Motiva                    & Port Arthur      & Jefferson    & 0.635  & 80  & 0    & 20  & 0   & 1648 \\
19 & Indorama Ventures Olefins & Port Neches      & Jefferson    & 0.2359 & 100 & 0    & 0   & 0   & 1648 \\
20 & LyondellBasell Industries & Corpus Christi   & Nueces       & 0.363  & 100 & 0    & 0   & 0   & 406  \\
21 & LyondellBasell Industries & Corpus Christi   & Nueces       & 0.767  & 15  & 30   & 35  & 20  & 406  \\
22 & MarkWest Javelina         & Corpus Christi   & Nueces       & 1      & 100 & 0    & 0   & 0   & 539  \\
23 & Dow Chemical Co.          & Orange           & Orange       & 0.791  & 100 & 0    & 0   & 0   & 1648 \\
24 & Exxon/SABIC               & Portland         & San Patricio & 1.8    & 100 & 0    & 0   & 0   & 516  \\
25 & Occidental Chem./Mexichem & Ingleside        & San Patricio & 0.544  & 100 & 0    & 0   & 0   & 539  \\
26 & Occidental Chem./Mexichem & Ingleside        & San Patricio & 0.773  & 15  & 30   & 35  & 20  & 539  \\
\bottomrule
\end{tabular}}
\caption{Location, ethylene production capacity, feedstock compositions of 26 ethylene plants in operation in Texas and their associated buses in ACTIVSg2000 dataset \cite{Birchfield2017}.} \label{tab:plants}
\end{table}

\begin{figure}[ht!]
    \centering
    \includegraphics[width=\textwidth]{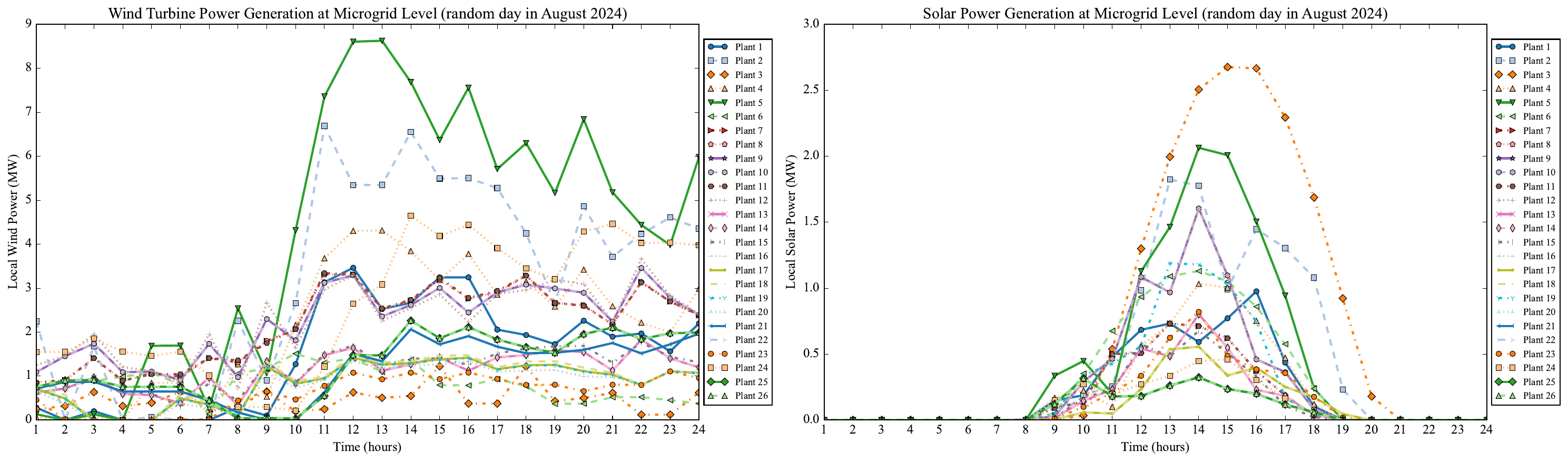}
    \vspace{-2em}
    \caption{Local hourly wind and solar availability profile for each ethylene plant microgrid on an arbitrary day in August 2024 (specifically, August 2, 2024).}
    \label{AugustREN}
\end{figure}

Recall that the microgrid associated with an ethylene plant can take VRE (solar and wind energy) locally. To accurately estimate the VRE power availability for each microgrid, we use the longitude and latitude values and find out the relevant local weather data, including ambient pressure and temperature, wind speed, and solar radiation intensity, of a given day from public or private weather API providers (e.g., visualcrossing dataset \cite{visualcrossing}). We then adopt the following relations to convert the weather data into the plant-level solar and wind power availability, denoted as $P^{\mathrm{WT}}$ and $P^{\mathrm{PV}}$, respectively:
\begin{equation}\label{eq:renewablecalculation}
    P^{\mathrm{WT}}_{p,t} = \frac{1}{2}\eta^{\mathrm{WT}}\pi r_{\mathrm{swept},p}^2\rho^{\mathrm{Air}}_{p,t}v^3_{\mathrm{Wind}, p,t}, \qquad P^{\mathrm{PV}}_{p,t} = A_{\mathrm{Panel},p} R^{\mathrm{Solar}}_{p,t},
\end{equation}
where the power coefficient $\eta^{\mathrm{WT}}$ for modern utility-scale wind turbine is around 50\% \cite{turbineefficiency}, $r_{\mathrm{swept}}$ denotes the swept radius of the locally installed wind turbine, $v_{\mathrm{Wind}}$ is the local wind speed (in m/s), $A_{\mathrm{Panel}}$ is total solar panel area (in \ce{m^2}) of the local solar farm, and $R^{\mathrm{Solar}}$ is the local solar radiation intensity (in \ce{MW/m^2}). The air density $\rho^{\mathrm{Air}}_{p,t}$ (in \ce{kg/m^3}) is determined using the ideal gas law as $\rho^{\mathrm{Air}}_{p,t} = \frac{P^{\mathrm{Air}}_{p,t}M^{\mathrm{Air}}}{RT_{p,t}}$, in which $M^{\mathrm{Air}}$ is the average molecular weight of air (28.97 \ce{g/mol}), $R$ is the ideal gas constant (8.314 \ce{J/K\cdot mol}), $P^{\mathrm{Air}}$ is the ambient pressure (in \ce{Pa}), and $T^{\mathrm{Air}}$ is the ambient temperature (in K). Using this approach, we can synthesize an hourly VRE availability profile for each of the 26 microgrids for any given day, as illustrated in Figure \ref{AugustREN}. These profiles will be sent to the joint optimization formulation as inputs. The specifications of fossil-based dispatchable units, battery storage, and hydrogen storage units present in the microgrid are adopted from \cite{Naraghi2025,khodaei2013microgrid,escape_25_microgrid} and are summarized in Table \ref{table_parameters}.

\begin{table}[ht!]
    \centering
    \begin{adjustbox}{width=\textwidth}
    \begin{tabular}{cccccc}
    \toprule
    Fuel group & \# of & Operating cost & Min-max & Min up/down & Ramp-up/down\\
    & units & (\$/MWh) & capacity (MW) & time (hr) &  rate (MW/h) \\
    \midrule
    Natural gas & 20 & 33.4 & 1-5 & 3 & 2.5 \\
    Hydrogen & 1 & 30 & $10^{-5}$-1 & N/A & N/A\\
    \bottomrule
    \toprule
    Storage & \# of & Storage cost & Capacity & Min-max charging/ & Min charging/ \\
    type & units & (\$/ton) & & discharging & discharging\\
    & & & & power (MW) & time (hr)\\
    \midrule
    Battery storage & 1 & N/A & 20 (MWh) & 0.8-4 & 5 \\
    Hydrogen storage & 1 & 10,000 & 10 (ton) & N/A & N/A \\
    \bottomrule
    \end{tabular}
    \end{adjustbox}
    \caption{Characteristics of fossil-based dispatchable units, battery storage, and hydrogen storage units present in the microgrid of Figure \ref{fig_superstructure}.} \label{table_parameters}
\end{table}

Now, in terms of the power system, we adopt the data sourced from the ACTIVSg2000 test case \cite{Birchfield2017}, a synthetic yet realistic power system benchmark. This test case includes 2000 buses and serves as an effective proxy for real-world large-scale systems, making it suitable for studying scalability and computational efficiency. Synthetic electric grid models are fictitious representations designed to be statistically and functionally similar to actual electric grids while containing no confidential critical energy infrastructure information. All 26 ethylene plant microgrids were mapped to buses based on geographic proximity, using their latitude and longitude data to determine their nearest grid connection points. These results are shown in Table \ref{tab:plants}. Note that some microgrids are connected to the same bus. In our problem, all 26 microgrids are connected to 15 distinct buses, leading to 15 LP-relexed subproblems in the Benders' decomposition procedure within our proposed two-stage solution methodology.


Finally, we recognize that different chemical plants typically adopt different strategies and pathways toward electrification. On one hand, large corporations may adopt a more progressive approach toward electrification, as they can invest more on R\&D activities and new technologies. Meanwhile, smaller companies may choose to be more conservative. On the other hand, large-scale plants may also be more reluctant to revamp their existing facilities with new ones due to the large capital investment needed compared to smaller plants. To consider different levels of electrification effort adopted by different ethylene plants without making the problem overwhelmingly complex, we categorize all 26 plants into three tiers based on their sizes. The first tier (T1), which represents the largest corporations, has only two ethylene plants whose ethylene production exceeds 3 million tons/year. The second tier (T2) contains 10 plants with ethylene production capacity falling in between 1 and 2.5 million tons/year. The third tier (T3), which produces less than 1 million tons of ethylene per year, consists of the remaining 14 plants. T1, T2, and T3 plants consititutes 25.3\%, 47.1\%, and 27.6\% of the total ethylene production capacity surveyed, respectively. Then, we vary the electrification level (i.e., percentage of ethylene produced from electrified crackers) for different tiers and formulate 10 cases shown in Table \ref{tab:case_study}. Case (0) represents the baseline case, whereas Cases (1) through (9) consider both homogeneous and heterogeneous electrification levels across all 26 microgrids.

\begin{table}[ht!]
    \centering
    \begin{tabular}{c|ccc}
        \toprule
        \multirow{2}{*}{\textbf{Case \#}} & \multicolumn{3}{c}{\textbf{Level of Electrification}} \\
        \cline{2-4}
        & \textbf{Tier 1} & \textbf{Tier 2} & \textbf{Tier 3} \\
        \midrule
        0 & 0\%  & 0\%  & 0\%  \\
        1 & 10\% & 10\% & 10\% \\
        2 & 20\% & 20\% & 20\% \\
        3 & 30\% & 30\% & 30\% \\
        4 & 40\% & 40\% & 40\% \\
        5 & 50\% & 50\% & 50\% \\
        6 & 50\% & 30\% & 10\% \\
        7 & 10\% & 30\% & 50\% \\
        8 & 10\% & 50\% & 50\% \\
        9 & 50\% & 10\% & 10\% \\
        \bottomrule
    \end{tabular}
    \caption{We formulate 10 cases by varying the electrification level for each tier. Also, for local solar power availability, we specify the solar panel area to be 8000 \ce{m^2} for T1 plants, 4000 \ce{m^2} for T2 plants, and 2000 \ce{m^2} for T1 plants. For local wind power availability, we specify the swept radius of wind turbines to be 100 \ce{m} for T1 plants, 50 \ce{m} for T2 plants, and 25 \ce{m} for T1 plants.}
    \label{tab:case_study}
\end{table}

\subsection{Computational Efficiency}

To evaluate the algorithmic performance of our proposed two-stage solution approach, we solve Cases (1) through (9) using the direct-solve method and the two-stage solution approach and compare the CPU time required to reach 0.1\% optimality gap. Note that different case studies only differ in input parameters, and all the cases will have the same size, with 238,034 continuous variables and 50,498 binary variables. Thus, this comparison directly reflects the differences in algorithm performance rather than differences in model size. All experiments are conducted using Gurobi 12.0.1 installed on a Dell Precision 7920 workstation equipped with Intel Xeon Gold 6226R 2.9 GHz CPU and 12$\times$8 GB DDR4 2933MHz RDIMM ECC memory. There are a couple of observations one can draw from the CPU time results in Table \ref{tab:time}. First, our two-stage method speeds up the optimization process in 7 out of 9 cases, achieving an averaging CPU time speedup of 93.5\% that is quite significant. Second, our two-stage method typically solves each case within 400-500 seconds, with the lowest CPU time being slightly less than 300 seconds (Case (2)) and the highest being slightly more than 600 seconds (Case (9)). Meanwhile, the direct-solve approach experiences much greater variations in CPU time, ranging from less than 200 seconds (Case (1)) to more than 20,000 seconds (Case (2)). This results in a much higher standard deviation in CPU time across different cases. In other words, our two-stage approach increases the computational time for easy-to-solve cases but greatly decreases the computational time for difficult-to-solve cases, thus homogenizing the computational efficiency by damping the case-wise variations. Overall, this study demonstrates that the proposed method offers considerable improvements in computational efficiency. 
 
\begin{table}[ht!]
    \centering
    \begin{tabular}{c|c|c|c}
        \toprule
        \multirow{2}{*}{\textbf{Case \#}} & \multicolumn{2}{c|}{\textbf{CPU Time (seconds)}} & \multirow{2}{*}{\textbf{\% Difference}} \\
        \cline{2-3}
        & \textbf{Direct-solve} & \textbf{Two-stage Method} & \\
        \hline
        1 & 197.6   & 401.2 & 103.0 \\
        2 & 21021.3 & 409.6 & -98.1 \\
        3 & 1028.6  & 264.1 & -74.3 \\
        4 & 16788.3 & 418.5 & -97.5 \\
        5 & 13079.8 & 479.0 & -96.3 \\
        6 & 1424.8  & 559.9 & -60.7 \\
        7 & 142.1   & 400.5 & 181.8 \\
        8 & 8136.7  & 475.4 & -94.2 \\
        9 & 784.7   & 633.1 & -19.3 \\
        \hline
        \bf Average & \bf 6956.0 & \bf 449.0 & \bf -93.5 \\
        \hline
        \bf Standard deviation & \bf 8134.2 & \bf 105.5 & \bf -98.7 \\
        \bottomrule
    \end{tabular}
    \caption{CPU time comparison between the direct-solve MILP model and the proposed two-stage solution methodology. The cases are defined in Table \ref{tab:case_study}.}
    \label{tab:time}
\end{table}

\subsection{Impacts of Electrification on Operating Costs and Emissions}


Without electrification (i.e., Case (0)), the minimum total operating cost of power system and microgrids is \$23.470 million per day. The operating cost breakdown is \$21.939 million per day (93.4\%) for the power system and \$1.531 million per day (6.5\%) for all 26 ethylene microgrids. In terms of environmental impacts, the direct or Scope 1 emissions, which are greenhouse gases released directly into the atmosphere from sources owned or controlled by an organization \cite{emissions}, are 1.2999 million tons per day for the base case, among which 91.9\% (or 1.1946 million tons/day) is associated with the power system sector and 8.1\% (or 0.1053 million ton/day) is associated with the ethylene plants. Clearly, the power system contributes the majority of operating costs and emissions compared to all ethylene plants combined, suggesting that the joint optimization that considers both the power systems and chemical plants is needed.

\subsubsection{Impacts on Operating Costs}

As we consider different levels of electrification, the optimal operating costs breakdown is summarized in Table \ref{tab:costs}. As the electrification level increases homogeneously from 10\% to 50\% across all 26 microgrids, the operating costs of both the power systems and microgrids increases, though at different rates. Specifically, despite having an order of magnitude larger scale in magnitude, the power systems operating costs increase only by \$0.9 million/day, whereas microgrids operating costs increase by \$5.73 million/day. Overall, the microgrid side contributes 86\% of the operating costs increase. This is expected as electrification occurs at the microgrid side in this study.

\begin{table}[ht!]
    \centering
    \begin{tabular}{c|c|c|c|c|c}
    \toprule
    \textbf{Case} & \textbf{Power System} & \textbf{Microgrids} & \textbf{Power Systems} & \textbf{Microgrids} & \textbf{Overall Costs} \\
    \bf \# & \bf  (million \$) & \bf  (million \$) & \bf \% Change & \bf \% Change & \bf \% Change \\
    \hline
    1 & 21.9146 & 2.1026 & -0.11 & 37.29 & 2.33 \\
    2 & 21.9146 & 2.7559 & -0.11 & 79.95 & 5.12 \\
    3 & 21.9146 & 3.4669 & -0.11 & 191.67 & 8.14 \\
    4 & 22.3480 & 5.2911 & 1.87 & 245.49 & 17.76 \\
    5 & 22.8397 & 7.2583 & 4.11 & 373.94 & 28.24 \\
    6 & 22.1519 & 4.1068 & 0.97 & 168.16 & 11.88 \\
    7 & 22.1566 & 4.1586 & 0.99 & 171.54 & 12.12 \\
    8 & 22.6056 & 5.9534 & 3.04 & 288.73 & 21.68 \\
    9 & 22.1848 & 3.3666 & 1.12 & 119.83 & 8.87 \\
    \bottomrule
    \end{tabular}
    \caption{Operational costs and percentage change rates for each case study relative to the baseline (Case~0).}
    \label{tab:costs}
\end{table}

In addition, we observe that there is be a sudden jump in optimal total operating costs after 30\% electrification. From 0\% to 30\% electrification, the optimal total operating costs rises gradually by less than 3\% for every 10\% increase in electrification level, whereas from 30\% to 50\% electrification, the optimal total operating costs increases by almost 9\% for every 10\% increase in electrification level. Specifically, between 10\% to 30\% electrification, the power systems operating costs are actually slightly lower than the base case and stay flat, whereas starting from 30\% electrification, both power systems and microgrids experience an increase in operating costs. The combined system of power systems and microgrids is much more expensive to operate once electrification is pushed beyond 30\%, primarily because the power systems side needs to commit more generators for more electrified microgrids, as local electricity generation becomes less sufficient.

Pairwise correlations between total operating costs and electrification levels are calculated to be 0.482 for Tier 1 plants, 0.924 for Tier 2 plants, and 0.811 for Tier 3 plants. Note that Tier 1 contains only the largest two plants, Tier 2 contains 10 medium-sized plants, and Tier 3 contains 14 smaller plants. Therefore, although Tier 2 dominates the aggregate cost response as it represents close to 50\% of the total ethylene production capacity, each Tier 1 plant has the greatest per-plant cost leverage. Furthermore, the Tier 3 aggregate effect is only slightly smaller than Tier 2 despite occupying a much smaller contribution to the overall ethylene production (27.6\% for Tier 3 vs.~47.1\% for Tier 2). This is likely because smaller plants often have weaker economies of scale, especially as the startup and shutdown costs of local generators are independent of plant size. In other words, electrifying a small plant may impose a larger operating cost penalty per ton of ethylene produced, making Tier 3 category look disproportionately important relative to its capacity share.

Another interesting observation is that the distribution of electrification level across different tiers have heterogeneous impact on optimal operating costs. For example, Cases (6) and (7) have capacity-weighted electrification levels near 30\% (29.5\% and 30.5\%, respectively). However, their optimal total operating costs are about \$0.9 million/day higher than the most similar homogeneous benchmark (Case (3)). Similarly, Case (8) has capacity-weighted electrification level of approximately 40\%, which is nearly the same as Case (4), but its optimal total operating cost is \$0.92 million/day higher. In summary, given similar capacity-weighted electrification level, we observe that heterogeneity in electrification efforts across different plants leads to more expensive operation compared to more balanced electrification efforts.

\subsubsection{Impacts on Emissions}

For each case, when the optimal solution is achieved, the associated system-wide daily Scopes 1 and 2 emissions and their breakdowns are summarized Tables \ref{tab:emissions} and \ref{tab:category}. Between 0 to 30\% electrification level, the system-wide emissions decrease as electrification level increases, and the largest system-wide Scope 1 emissions reduction (1.357\%) occurs at 30\% electrification for all tiers (Case (3)). While this percentage reduction might seems small, it is actually quite significant, considering that an annual \ce{CO_2}-equivalent Scope 1 emissions reduction of 6.424 million tons can be achieved by just electrifying steam cracking process alone. 

\begin{table}[ht!]
    \centering
    \begin{tabular}{c|c|c}
    \toprule
       \textbf{Case \#} & \textbf{Scope 1 Emissions} & \textbf{\% Difference Compared}  \\
       & \textbf{(million tons/day)} & \textbf{to the Base Case}\\
       \hline
        0 & 1.2999 & N/A    \\
        \hline
        1 & 1.2946 & -0.404 \\
        2 & 1.2883 & -0.890 \\
        3 & 1.2823 & -1.357 \\
        4 & 1.2933 & -0.510 \\
        5 & 1.3050 & 0.391  \\
        6 & 1.2931 & -0.521 \\
        7 & 1.2913 & -0.564 \\
        8 & 1.3025 & 0.197  \\
        9 & 1.2984 & -0.118 \\
        \bottomrule
    \end{tabular}
    \caption{Total \ce{CO_2}-equivalent emissions released for each case study.}
    \label{tab:emissions}
\end{table}

\begin{table}[ht!]
    \centering
    \begin{adjustbox}{max width=\textwidth}
    \begin{tabular}{c|cc|cc}
        \toprule
        \bf Case & \multicolumn{2}{c}{\textbf{Scope 1 Emissions (Mtons/day)}} & \multicolumn{2}{c}{\textbf{Scope 2 Emissions (Mtons/day)}}\\
        \cline{2-5}
        \bf \# & \textbf{Power Systems} & \textbf{Microgrids} & \textbf{Ethylene Production} & \textbf{Other Demand} \\
        \hline
        0 & 1.1947 & 0.1053 & 0.1053 & 1.1947 \\
        1 & 1.1955 & 0.0991 & 0.1000 & 1.1946 \\
        2 & 1.1962 & 0.0922 & 0.0937 & 1.1947 \\
        3 & 1.1962 & 0.0861 & 0.0876 & 1.1946 \\
        4 & 1.2165 & 0.0767 & 0.0971 & 1.1962 \\
        5 & 1.2393 & 0.0657 & 0.1081 & 1.1969 \\
        6 & 1.2084 & 0.0847 & 0.0974 & 1.1958 \\
        7 & 1.2074 & 0.0839 & 0.0955 & 1.1958 \\
        8 & 1.2282 & 0.0742 & 0.1063 & 1.1962 \\
        9 & 1.2080 & 0.0904 & 0.1026 & 1.1958 \\
        \bottomrule
    \end{tabular}
    \end{adjustbox}
    \caption{Scope 1 and 2 \ce{CO2}-equivalent emissions for power systems and microgrid stakeholders. Scope 2 emissions are indirect GHG emissions associated with the purchase of electricity, steam, heat, or cooling \cite{emissions}. Scope 2 emissions can be categorized based on the end use: emissions associated with ethylene production and those from other demand met by the power system. We remark that, by definition, the emissions resulting from the power transferred from main grid to the microgrids are not included in the microgrid emissions. Thus, the Scope 2 emissions associated with ethylene production will always be higher than the corresponding Scope 1 emissions associated with the microgrids.}
    \label{tab:category}
\end{table}

Meanwhile, Table \ref{tab:emissions} reveals another important insight; that is, higher levels of electrification do not always lead to lower emissions. In this case, beyond 30\% electrification, the system-wide Scopes 1 and 2 emissions both increase. This is consistent with the observation we made in our earlier work \cite{Naraghi2025,escape_25_microgrid}. Specifically, between 0\% to 50\% electrification, the average natural gas consumption by conventional crackers declines as expected. Nevertheless, the rate of decrease is the greatest between 0\% to 30\% electrification and then diminishes beyond 30\% electrification. On the other hand, at the microgrid level, the natural gas consumed by local dispatchable generators shows a steady increase as the level of electrification increases. Meanwhile, when renewable electricity generation technologies (e.g., fuel cells) are not yet mature enough to be cost-competitive compared to conventional fossil-based generators, as the power demand from electrified cracker units increases, the microgrid must produce more electricity using its local generators or withdraw more power from the grid (which is also less renewable). In the former case, due to efficiency losses, this will incur more energy demand than directly using natural gas as fuel for conventional crackers, thereby incurring more costs and higher carbon intensity. Overall, we expect that a higher electrification level may become favorable from a cost and emission perspective only when clean energy sources and technologies become more cost effective to acquire and operate.

\begin{figure}[ht!]
    \centering

    \begin{minipage}[t]{0.32\textwidth}
        \centering
        \includegraphics[width=\textwidth]{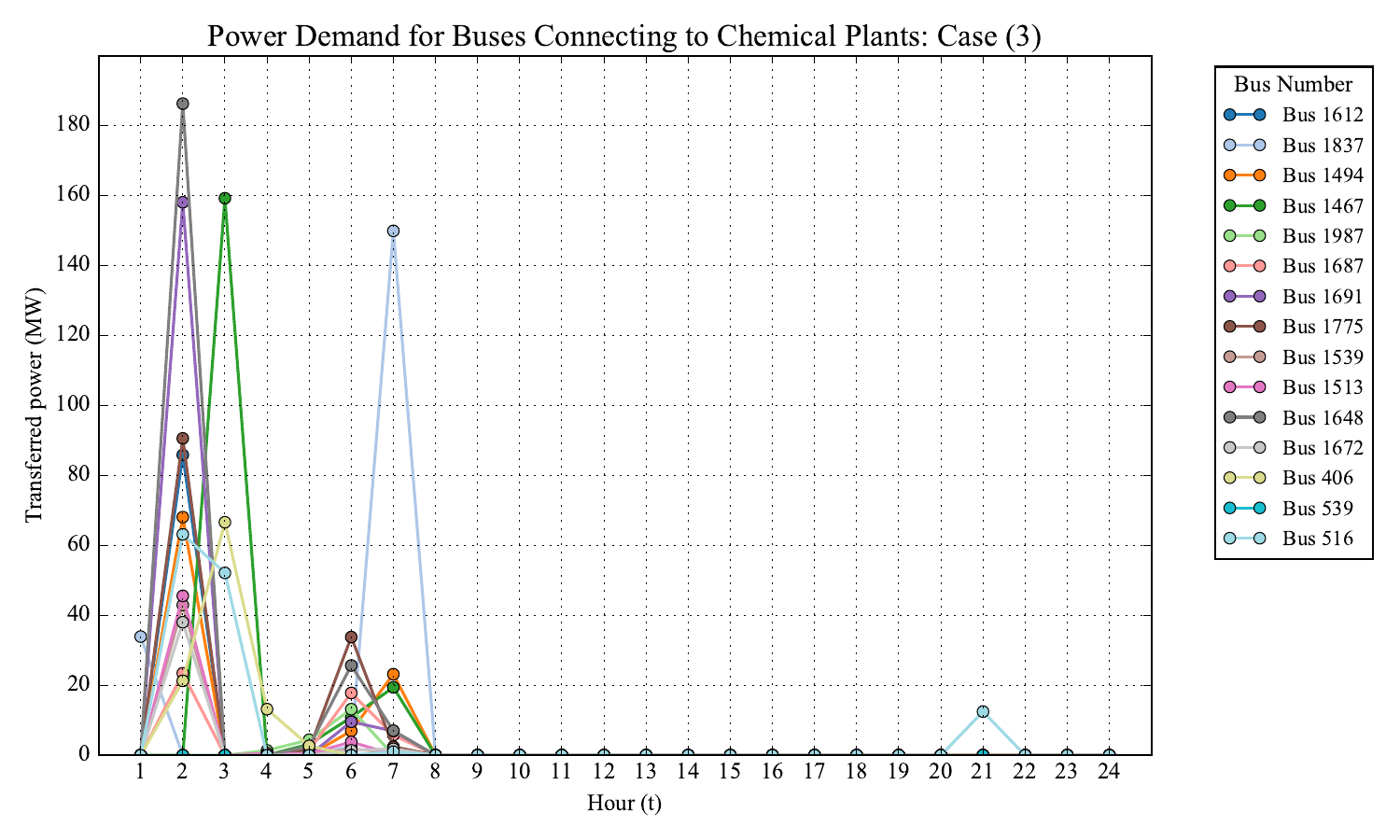}
    \end{minipage}
    \hfill
    \begin{minipage}[t]{0.32\textwidth}
        \centering
        \includegraphics[width=\textwidth]{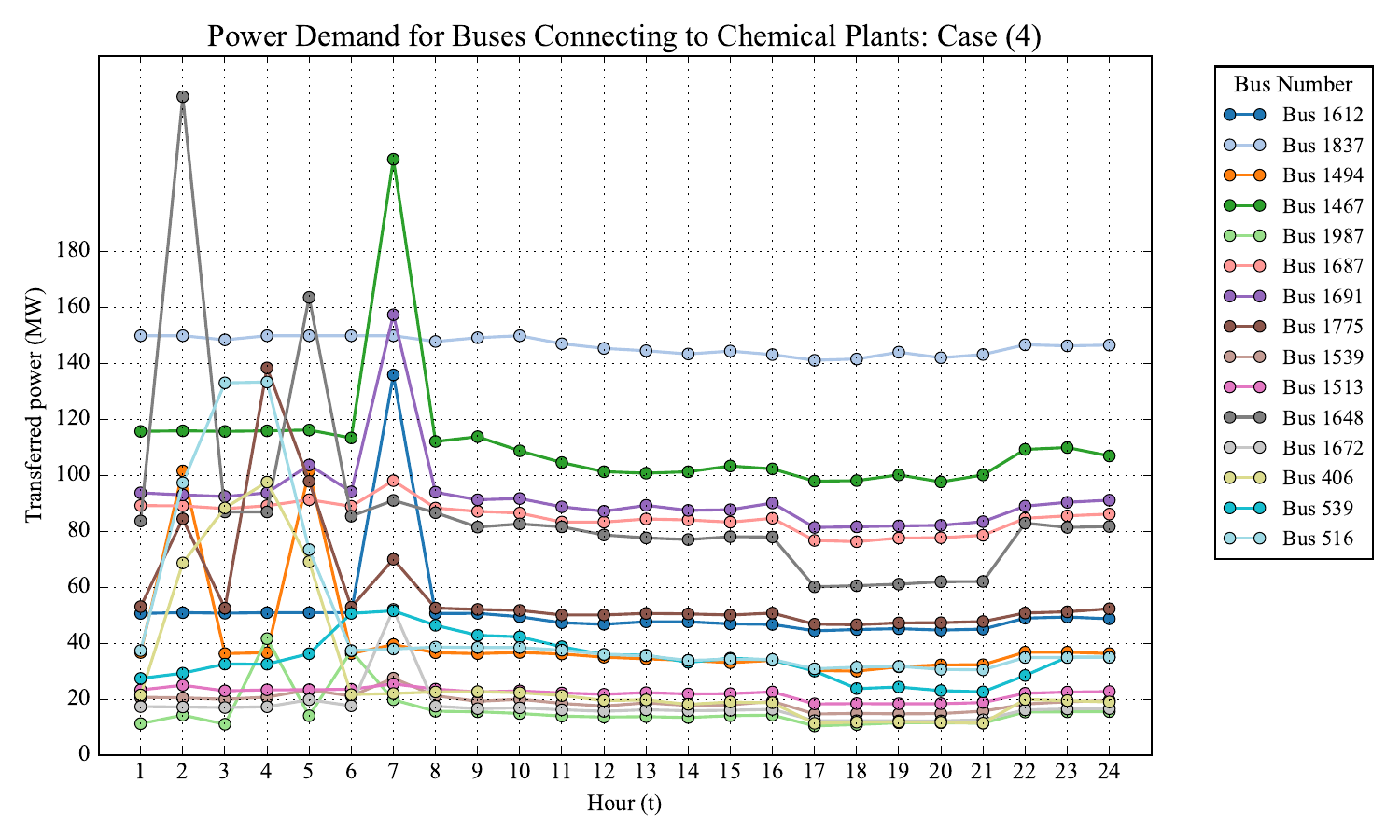} 
    \end{minipage}
    \hfill
    \begin{minipage}[t]{0.32\textwidth}
        \centering
        \includegraphics[width=\textwidth]{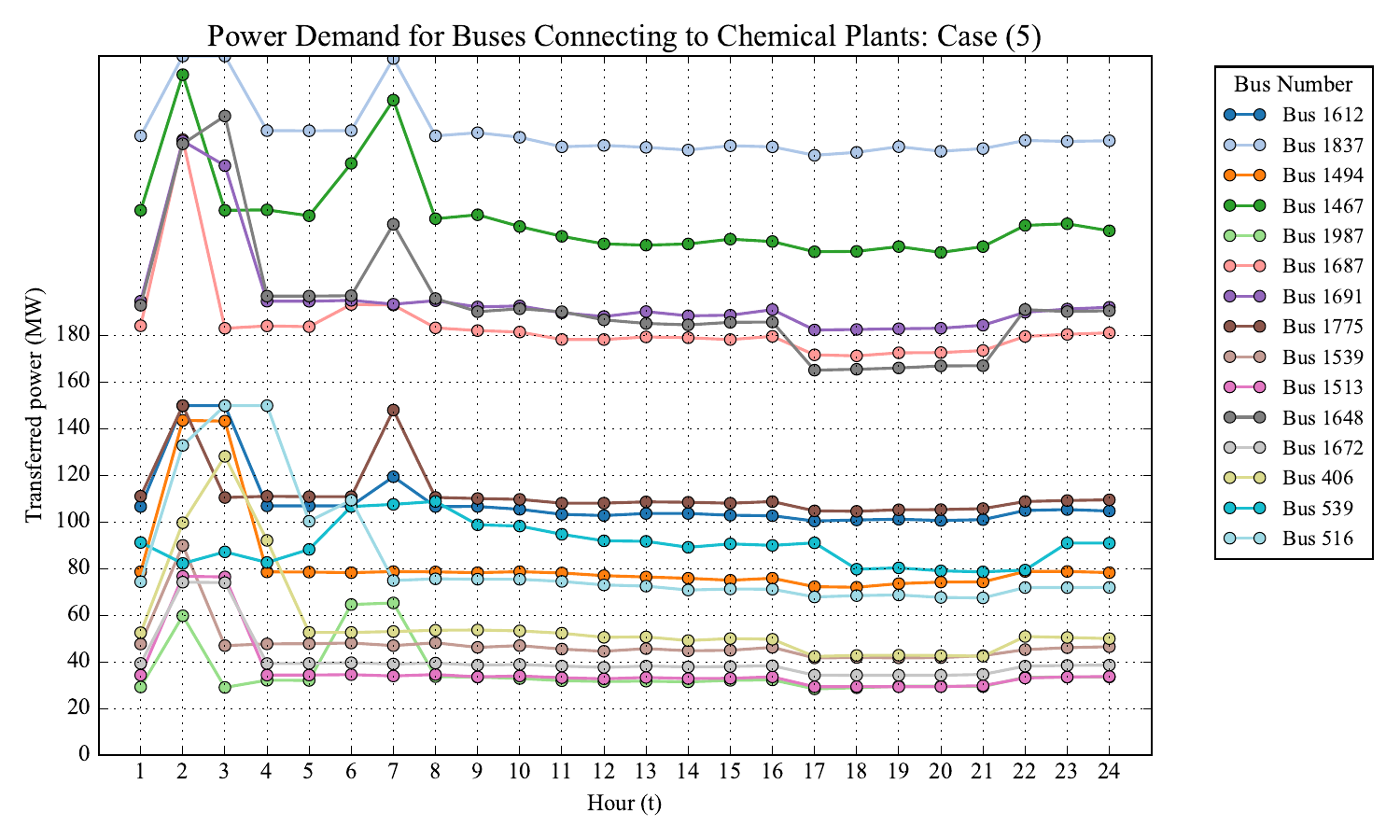} 
    \end{minipage}
    \vspace{-1em}
    \caption{Power transferred from power systems to microgrids for 30\%, 40\%, and 50\% electrification.}
    \label{Transferred_Power}
\end{figure}

Furthermore, Figure \ref{Transferred_Power} illustrates how the power system and microgrids interact over time. During non-peak hours when electricity market prices are low, more power is withdrawn from the main grid to the microgrids and stored. During peak hours with higher prices, microgrids prefer to generate electricity locally and/or use the electricity stored to meet the energy demand while minimizing operating costs. In particular, at 30\% electrification, most microgrids operate in islanded mode during peak hours. Meanwhile, starting at 40\% electrification, we observe that all microgrids operate at grid-connected mode during every hour of the day, as local fossil-based or renewable electricity generation becomes insufficient. 

\subsection{Spatial and Temporal Dependence}

\begin{figure}[ht!]
    \centering
    \includegraphics[width = 1\textwidth]{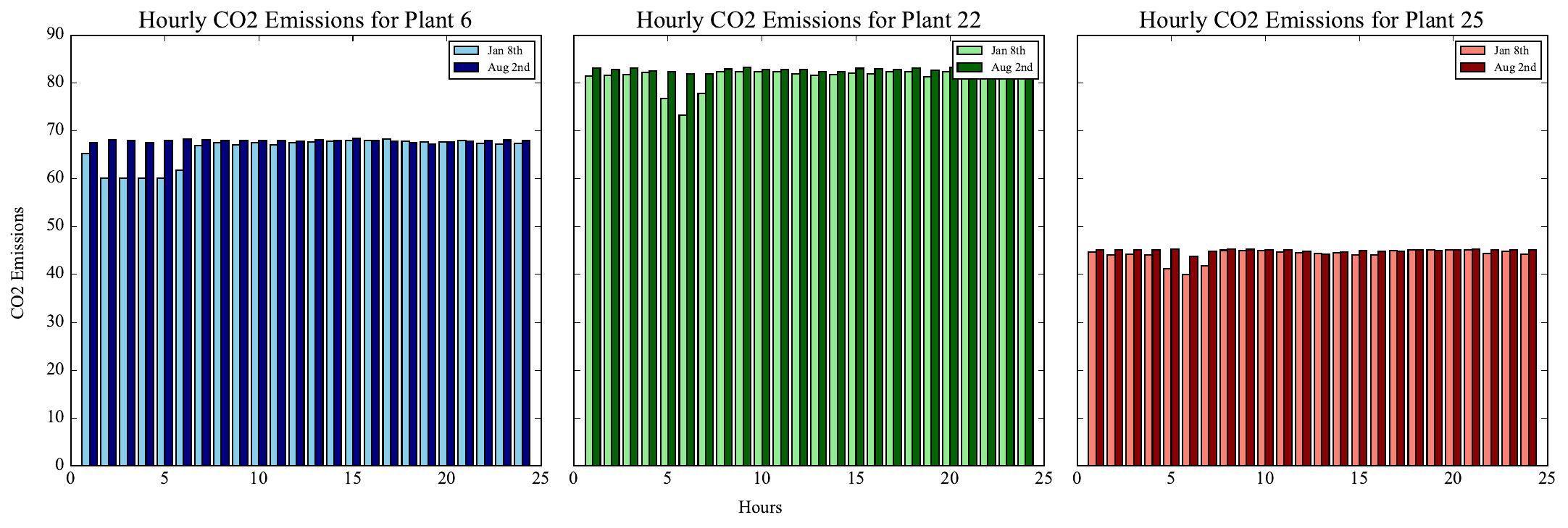}
    \vspace{-1em}
    \caption{Hourly Scope 1 \ce{CO_2}-equivalent emissions profiles for Plants 6, 22, and 25 on January 8 and August 2, 2024. }
    \label{PlantCO2}
\end{figure}

The VRE availability exhibits strong spatial and temporal dependence, as well as seasonality change. To shed light on how these variations affect the operation of microgrids, we select three representative microgrids, Plants 6, 22, and 25, and two representative days (January 8 for winter and August 2 for summer) in 2024 and examine their local hourly emissions profiles, which are shown in Figure \ref{PlantCO2}. All microgrids have 30\% electrification level (Case (3)). Plants 22 and 25 are less than 12 miles (19 km) from each other, producing 1 million tons and 0.544 million tons of ethylene each year, respectively. Meanwhile, Plant 6 has a similar ethylene production capacity (0.64 million ton/year) as Plant 25, but is more than 300 miles (482 km) away from Plants 22 or 25. 

\begin{figure}[ht!]
\centering
\begin{subfigure}{0.32\textwidth}
    \centering
    \includegraphics[width=\linewidth]{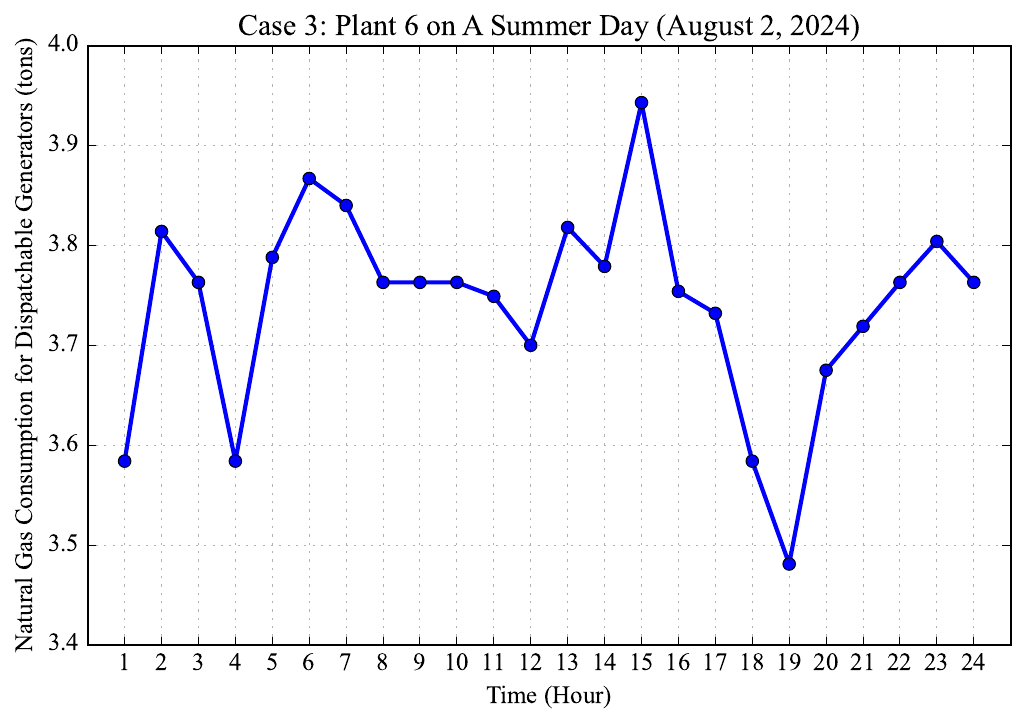}
\end{subfigure}
\hfill
\begin{subfigure}{0.32\textwidth}
    \centering
    \includegraphics[width=\linewidth]{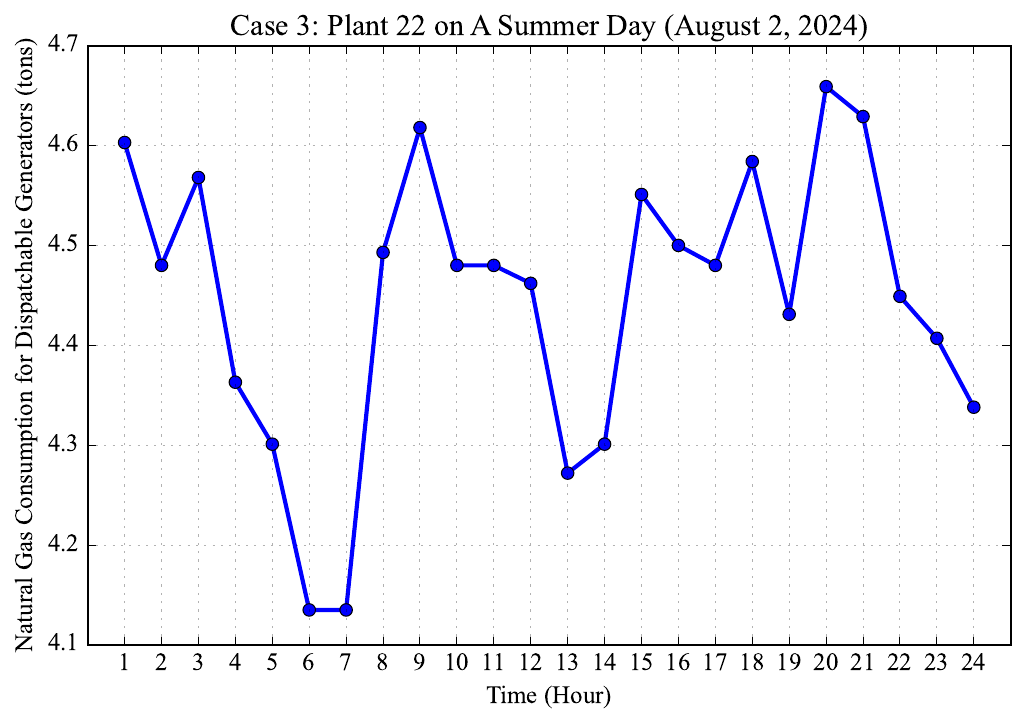}
\end{subfigure}
\hfill
\begin{subfigure}{0.32\textwidth}
    \centering
    \includegraphics[width=\linewidth]{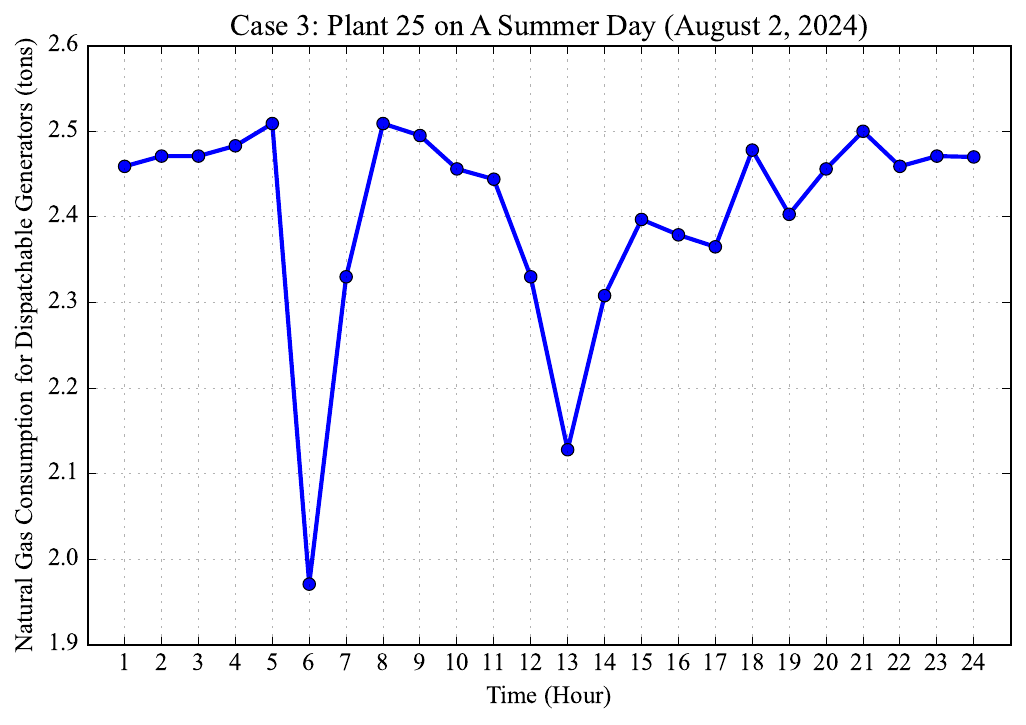}
\end{subfigure}
\bigskip
\begin{subfigure}{0.32\textwidth}
    \centering
    \includegraphics[width=\linewidth]{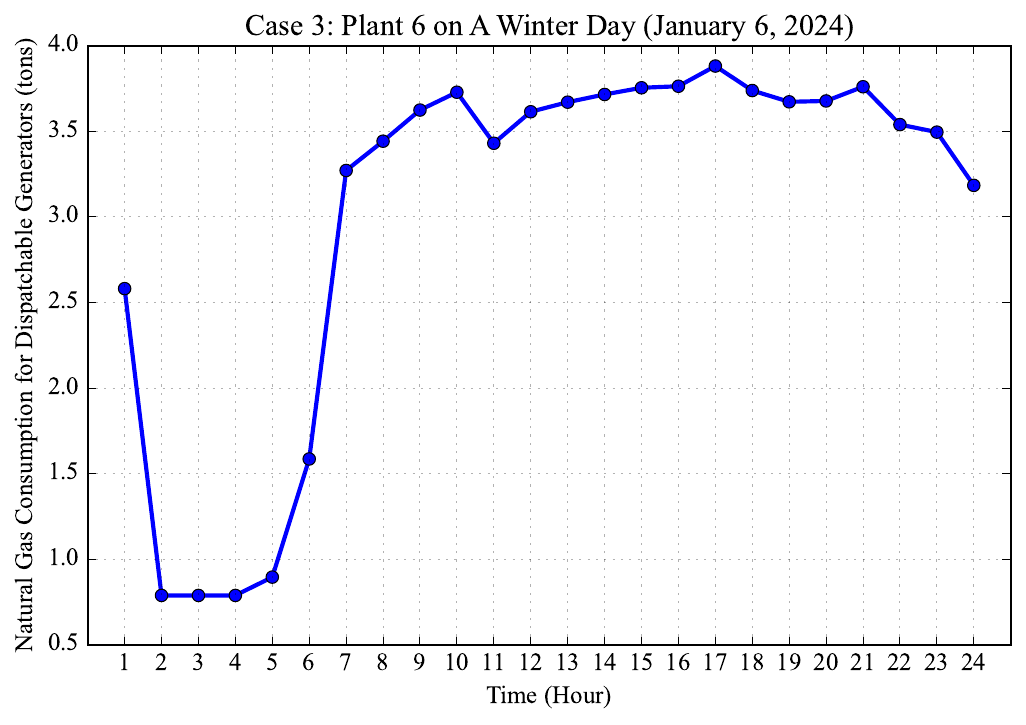}
\end{subfigure}
\hfill
\begin{subfigure}{0.32\textwidth}
    \centering
    \includegraphics[width=\linewidth]{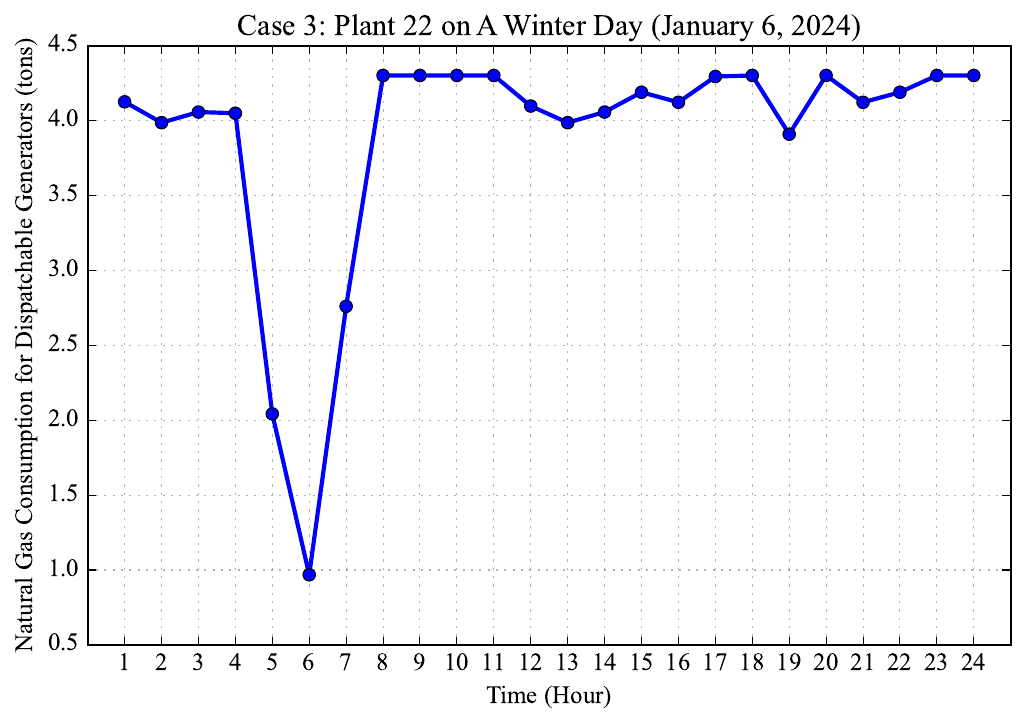}
\end{subfigure}
\hfill
\begin{subfigure}{0.32\textwidth}
    \centering
    \includegraphics[width=\linewidth]{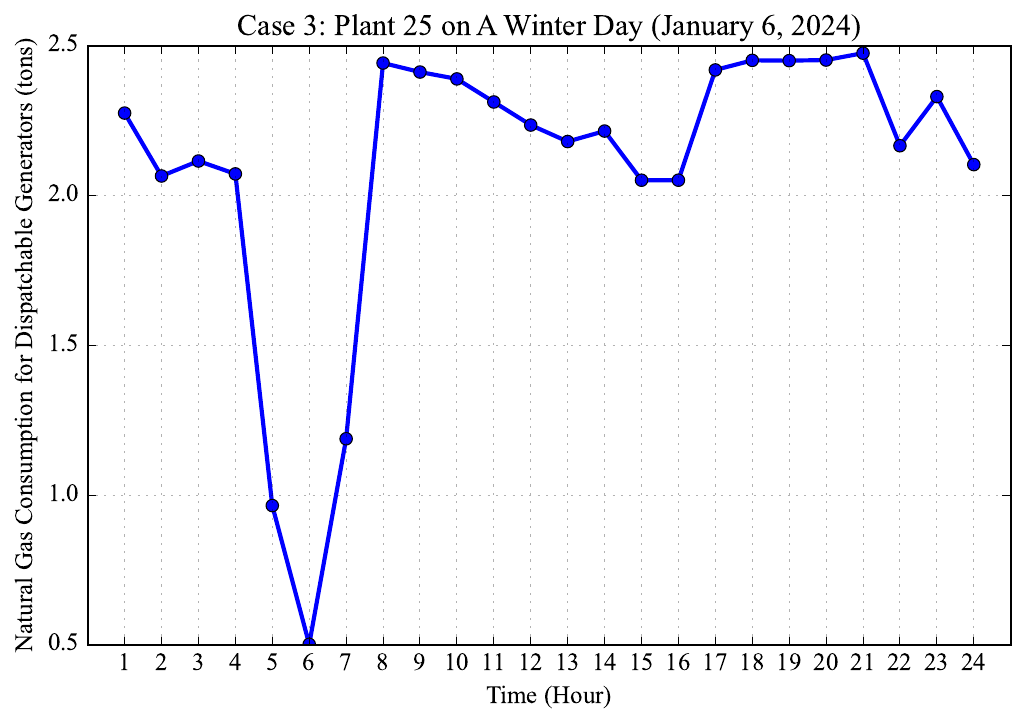}
\end{subfigure}
\caption{Natural gas consumption profile for local dispatchable generators in Plants 6 (left column), 22 (middle column), and 25 (right column) on a typical summer (top row) and winter day (bottom row). It is clear that plants nearby (Plants 22 and 25), despite having significant difference in production capacity, exhibit closer similarity in natural gas consumption profile compared to Plant 6 in both summer and winter days. Meanwhile, all three plants exhibit seasonal variations in local natural gas consumption profiles.}
\label{fig:plants_operation_NG}
\end{figure}

\begin{figure}[ht!]
\centering
\begin{subfigure}{0.32\textwidth}
    \centering
    \includegraphics[width=\linewidth]{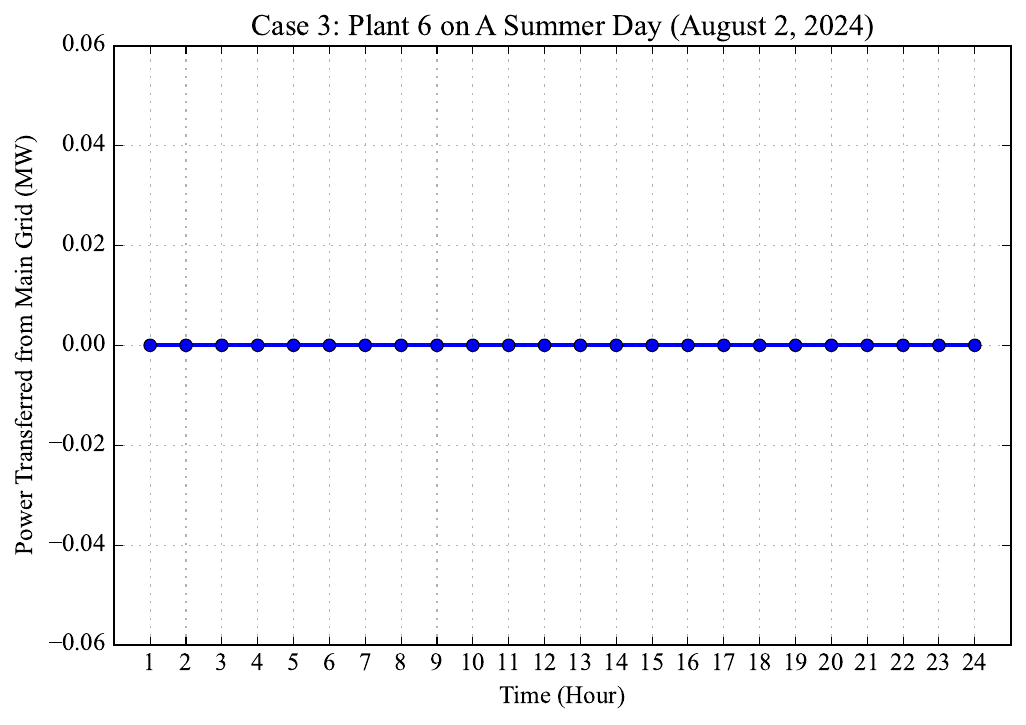}
\end{subfigure}
\hfill
\begin{subfigure}{0.32\textwidth}
    \centering
    \includegraphics[width=\linewidth]{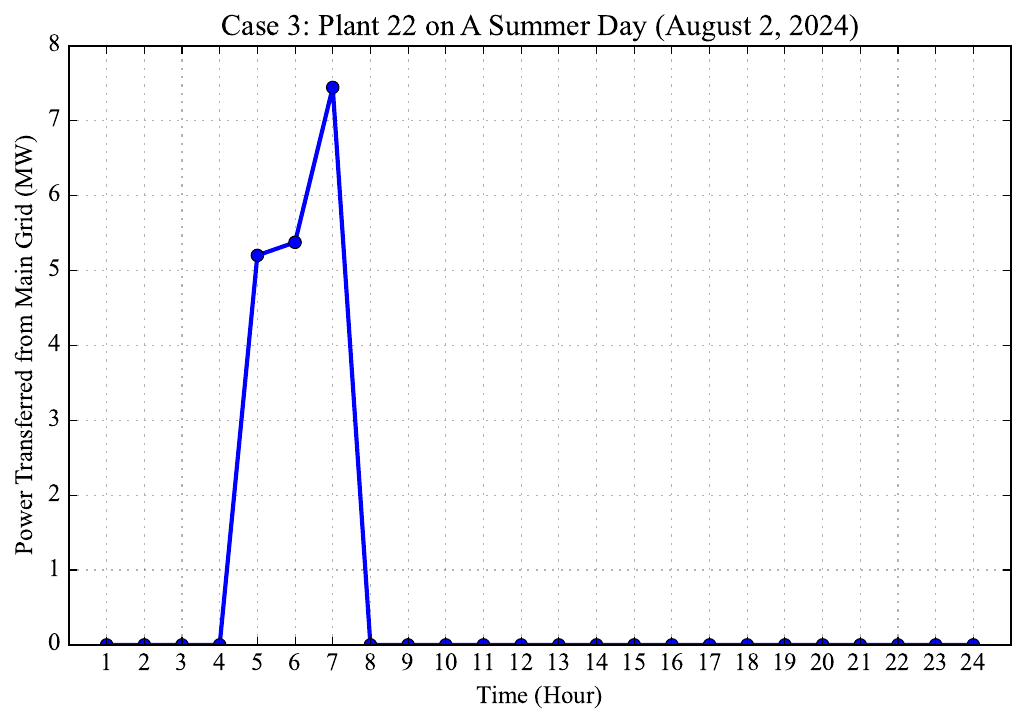}
\end{subfigure}
\hfill
\begin{subfigure}{0.32\textwidth}
    \centering
    \includegraphics[width=\linewidth]{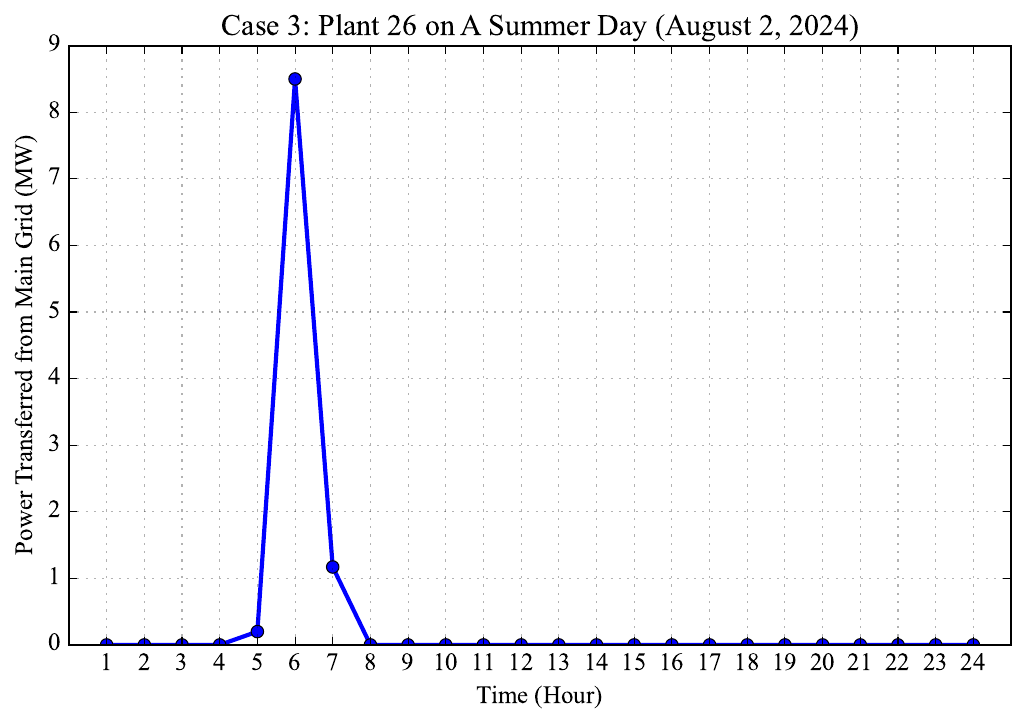}
\end{subfigure}
\bigskip
\begin{subfigure}{0.32\textwidth}
    \centering
    \includegraphics[width=\linewidth]{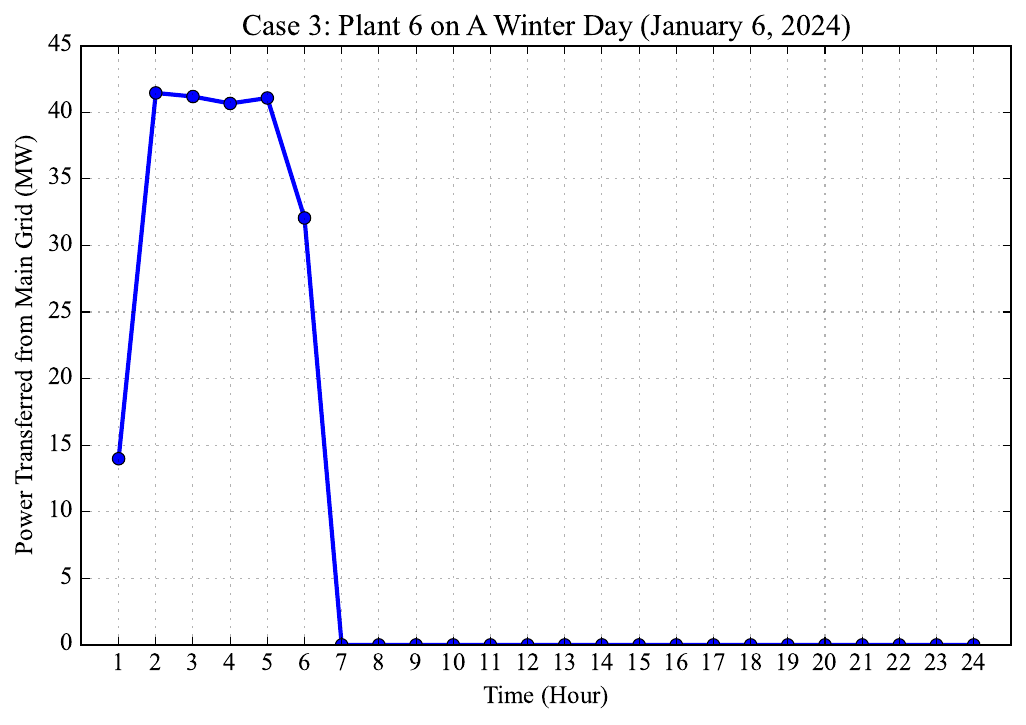}
\end{subfigure}
\hfill
\begin{subfigure}{0.32\textwidth}
    \centering
    \includegraphics[width=\linewidth]{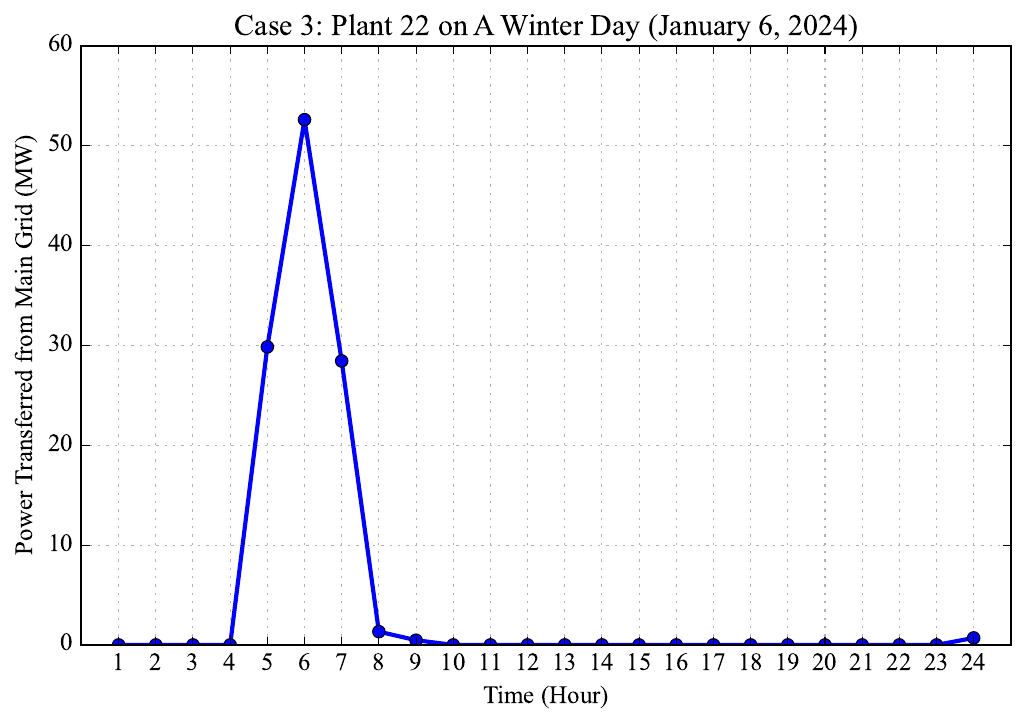}
\end{subfigure}
\hfill
\begin{subfigure}{0.32\textwidth}
    \centering
    \includegraphics[width=\linewidth]{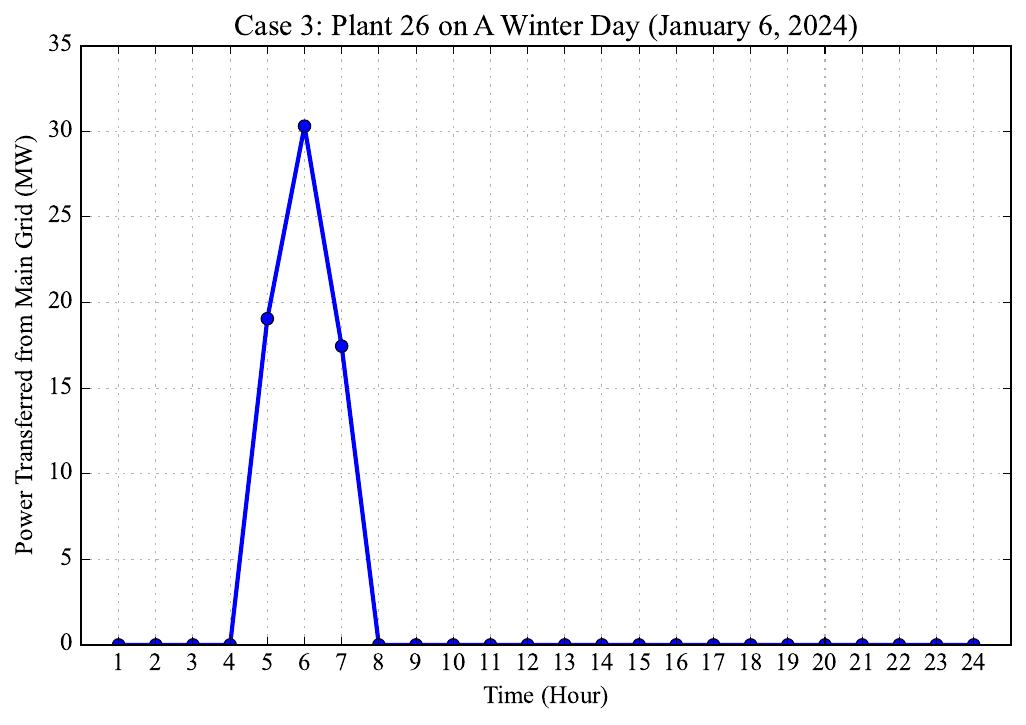}
\end{subfigure}
\caption{Power transferred from the main power grid to the microgrid for Plants 6 (left column), 22 (middle column), and 25 (right column) on a typical summer (top row) and winter day (bottom row). Like in Figure \ref{fig:plants_operation_NG}, we also observe strong spatial and temporal patterns in total power transferred from the main grid.}
\label{fig:plants_operation_Power}
\end{figure}

From Figures \ref{PlantCO2} through \ref{fig:plants_operation_Power}, it is clear that there is a strong temporal dependency in microgrid operation. Such temporal dependency is not only in daily scale but also in seasonal scale. Furthermore, since Plants 6 and 22 are both T3 plants, they are specified with the same solar and wind generation capacities. However, we observe that the hourly emissions profiles for Plants 6 and 22 are quite different, whereas those for Plants 22 and 25 are quite similar. This observation also matches with the power withdrawal and local dispatchable unit operational profiles shown in Figures \ref{fig:plants_operation_NG} and \ref{fig:plants_operation_Power}. In other words, plants located nearby have very similar operation and emission profiles even if they have very different production capacities. We remark that being able to capture and model these spatial and temporal dependencies is important because they provide deep insights on how extreme weather can affect the stability of local nodes in a power grid, which can significantly impact the power systems operation.

\section{Conclusion}

We present a multi-agent UCP framework for jointly optimizing power systems operation and electrified steam cracking microgrids. The proposed formulation captures the operational coupling between grid-level unit commitment and microgrid-level operation and energy management in a novel superstructure for steam cracking electrification. To improve tractability for large-scale instances, we propose a two-stage solution strategy, in which Benders decomposition with LP-relaxed microgrid subproblems provides a high-quality warm start for the full centralized MILP. Application to 26 ethylene plants in Texas connected to a 2000-bus synthetic grid network demonstrates that the proposed framework can solve realistic industrial electrification planning problems with substantially improved computational reliability relative to direct MILP solution.

The first realistic full-scale case study shows that electrification of steam cracking has non-monotonic economic and environmental impacts. Moderate electrification reduces system-wide emissions, with the largest Scope 1 emissions reduction occurring at 30\% electrification across all plant tiers. However, further electrification beyond this level increases both operating costs and emissions because local renewable generation, conversion, and storage are insufficient and expensive to meet the additional power demand without increased reliance on grid electricity and local fossil-based dispatchable generation. The increases in operating costs are dominated by the chemical plant microgrid side, and heterogeneous electrification across plant tiers is more expensive than balanced electrification at similar capacity-weighted electrification levels. These results indicate that electrification alone does not guarantee decarbonization. Instead, its benefits depend critically on the availability of clean, cost-effective electricity and coordinated operation between chemical plants and the power system.

The results also reveal important plant-size, spatial, and temporal dependencies. Tier 2 plants dominate the aggregate cost response because they represent the largest share of ethylene production capacity, while Tier 3 plants show disproportionately large cost impacts relative to their capacity share, likely due to weaker economies of scale and size-independent startup and shutdown costs. Local renewable availability and microgrid operation vary significantly across both location and season, leading to distinct emissions and dispatch profiles among geographically separated plants. Overall, this study highlights the need for coordinated, spatially resolved, and grid-aware electrification strategies for chemical manufacturing.

In our future research, we will focus on extending the framework in several directions. First, incorporating decentralized or privacy-preserving optimization methods would enable coordination among plants while safeguarding proprietary operational data, making the approach more applicable to real-world competitive markets. Second, integrating dynamic electricity market participation, real-time uncertainties, and investment decisions for low-carbon technologies would enhance the model's robustness and policy relevance. Together, these extensions would create a more robust and comprehensive decision support tool that better align industrial decarbonization strategies with cost-effectiveness, resilience, and privacy-preserving properties.

\section*{Acknowledgements}
This work is supported by the U.S. National Science Foundation under Award 2343072.

\printnomenclature
\bibliography{Ref}
\bibliographystyle{ieeetr}

\end{document}